\RequirePackage{ifpdf}
\ifpdf 
\documentclass[pdftex]{sigma}
\else
\documentclass{sigma}
\fi

\def\R{\mathbb{R}}               
\def\Z{\mathbb{Z}}      
\def\lcf{\lbrack\! \lbrack}
\def\rcf{\rbrack\! \rbrack}
\begin{document}

\allowdisplaybreaks

\renewcommand{\PaperNumber}{049}

\renewcommand{\thefootnote}{$\star$}

\FirstPageHeading

\ShortArticleName{Reduction of Symplectic Lie Algebroids}

\ArticleName{Reduction of Symplectic Lie Algebroids\\ by a Lie Subalgebroid and a Symmetry Lie Group\footnote{This paper is a contribution to the Proceedings
of the Workshop on Geometric Aspects of Integ\-rable Systems
 (July 17--19, 2006, University of Coimbra, Portugal).
The full collection is available at
\href{http://www.emis.de/journals/SIGMA/Coimbra2006.html}{http://www.emis.de/journals/SIGMA/Coimbra2006.html}}}

\Author{David IGLESIAS~$^\dag$, Juan Carlos MARRERO~$^\ddag$, David MART{\'I}N
DE DIEGO~$^\dag$,\\
 Eduardo MART{\'I}NEZ~$^\S$ and Edith PADR{\'O}N~$^\ddag$}

\AuthorNameForHeading{D. Iglesias, J.C. Marrero, D. Mart{\'{\i}}n de Diego, E. Mart{\'{\i}}nez and E. Padr{\'o}n}

\Address{$^\dag$~Departamento de Matem\'aticas, Instituto de Matem\'aticas y F{\'\i}sica Fundamental,\\
$\phantom{^\dag}$~Consejo Superior de Investigaciones Cient\'{\i}\-ficas, Serrano 123, 28006 Madrid, Spain}
 \EmailD{\href{mailto:iglesias@imaff.cfmac.csic.es}{iglesias@imaff.cfmac.csic.es}, \href{mailto:d.martin@imaff.cfmac.csic.es}{d.martin@imaff.cfmac.csic.es}}

\Address{$^\ddag$~Departamento de Matem\'atica Fundamental, Facultad de Ma\-te\-m\'a\-ti\-cas,\\
$\phantom{^\ddag}$~Universidad de la Laguna, La Laguna, Tenerife, Canary Islands, Spain}
\EmailD{\href{mailto:jcmarrer@ull.es}{jcmarrer@ull.es}, \href{mailto:mepadron@ull.es}{mepadron@ull.es}}

\Address{$^\S$~Departamento de Matem\'atica Aplicada, Facultad de Ciencias,
Universidad de Zaragoza,\\
$\phantom{^\S}$~50009 Zaragoza, Spain}
\EmailD{\href{mailto:emf@unizar.es}{emf@unizar.es}}

\ArticleDates{Received November 14, 2006, in f\/inal form March
06, 2007; Published online March 16, 2007}

\Abstract{We describe the  reduction procedure for a symplectic
Lie algebroid by a Lie sub\-algebroid and a symmetry Lie group.
Moreover, given an invariant Hamiltonian function we obtain the
corresponding reduced Hamiltonian dynamics. Several examples
illustrate the generality of the theory.}


\Classification{53D20; 58F05; 70H05}

\section{Introduction}
As it is well known, symplectic manifolds play a fundamental role
in the Hamiltonian formulation of Classical Mechanics. In fact, if
$M$ is the conf\/iguration space of a classical mechanical system,
then the phase space is the cotangent bundle $T^*M$, which is
endowed with a canonical symplectic structure, the main ingredient
to develop Hamiltonian Mechanics. Indeed, given a Hamiltonian
function $H$, the integral curves of the corresponding Hamiltonian
vector f\/ield are determined by the Hamilton equations.

A method to obtain new examples of symplectic manifolds comes from
dif\/ferent reduction procedures. One of these procedures is the
classical Cartan symplectic reduction process: \emph{If
$(M,\omega)$ is a symplectic manifold and $i:C\to M$ is a
coisotropic submanifold of $M$ such that its characteristic
foliation ${\mathcal F}=\ker (i^*\omega)$ is simple, then the
quotient manifold $\pi:C\to C/{\mathcal F}$ carries a unique
symplectic structure $\omega_r$ such that
$\pi^*(\omega_r)=i^*(\omega)$}.

A particular situation of it is the well-known Marsden--Weinstein
reduction in the presence of a $G$-equivariant momentum map
\cite{MW}. In fact, in \cite{Gotay-Tuynman} it has been proved
that, under mild assumptions, one can obtain any symplectic
manifold as a result of applying a Cartan reduction of the
canonical symplectic structure on $\R ^{2n}$. On the other hand,
an interesting application in Mechanics is the case when we have a
Hamiltonian function on the symplectic manifold satisfying some
natural conditions, since one can reduce the Hamiltonian vector
f\/ield to the reduced symplectic manifold and therefore, obtain a
reduced Hamiltonian dynamics.

A category which is closely related to symplectic and Poisson
manifolds is that of Lie algebroids. A Lie algebroid is a notion
which unif\/ies tangent bundles and Lie algebras, which suggests its
relation with Mechanics. In fact, there has been recently a lot of
interest in the geometric description of Hamiltonian (and
Lagrangian) Mechanics on Lie algebroids (see, for instance,
\cite{GGU,LMM,M,weinstein}). An important application of this
Hamiltonian Mechanics, which also comes from reduction, is the
following one: if we consider a principal $G$-bundle $\pi : Q \to
M$ then one can prove that the solutions of the
Hamilton-Poincar\'e equations for a $G$-invariant Hamiltonian
function $H : T^*Q \to \R$ are just the solutions of the Hamilton
equations for the reduced Hamiltonian $h : T^*Q/G \to \R$ on the
dual vector bundle $T^*Q/G$ of the Atiyah algebroid
$\tau_{TQ/G}:TQ/G\to \R$ (see \cite{LMM}).

Now, given a Lie algebroid $\tau_A:A\to M$, the role of the
tangent of the cotangent bundle of the conf\/iguration manifold is
played by $A$-tangent bundle to $A^*$, which is the subset of
$A\times TA^*$ given by
\begin{gather*}
{\mathcal T}^AA^*=\{(b,v)\in A\times
TA^*/\rho_A(b)=(T\tau_{A^*})(v)\},
\end{gather*}
where $\rho_A$ is the anchor map of $A$ and $\tau_{A^*}:A^*\to M$
is the vector bundle projection of the dual bundle $A^*$ to $A$.
In this case, ${\mathcal T}^AA^*$ is a Lie algebroid over $A^*$.
In fact, it is the  pullback  Lie algebroid of $A$ along the
bundle projection $\tau_{A^*}:A^*\to M$ in the sense of Higgins
and Mackenzie~\cite{HM}. Moreover, ${\mathcal T}^AA^*$ is a
symplectic Lie algebroid, that is, it is endowed with
a~nondegenerate and closed 2-section. Symplectic Lie algebroids are
a natural generalization of symplectic manifolds since, the f\/irst
example of a symplectic Lie algebroid is the tangent bundle of a
symplectic manifold.

The main purpose of this paper is to describe a reduction
procedure, analogous to Cartan reduction, for a symplectic Lie
algebroid in the presence of a Lie subalgebroid and a symmetry Lie
group. In addition, for Hamiltonian functions which satisfy some
invariance properties, it is described the process to obtain the
reduced Hamiltonian dynamics.

The paper is organized as follows. In Section 2, we recall the
def\/inition of a symplectic Lie algebroid and describe several
examples which will be useful along the paper. Then, in addition,
we describe how to obtain Hamilton equations for a symplectic Lie
algebroid and a Hamiltonian function on it.

\looseness=1
Now, consider a symplectic Lie algebroid $\tau _A:A\to M$ with
symplectic 2-section $\Omega_A$ and $\tau _B:B\to N$ a Lie
subalgebroid of $A$. Then, in Section 3 we obtain our main result.
Suppose that a Lie group $G$ acts properly and free on $B$ by
vector bundle automorphisms. Then, if $\Omega_B$ is the
restriction to $B$ of $\Omega_A$, we obtain conditions for which
the reduced vector bundle $\tau
_{\widetilde{B}}:\widetilde{B}=(B/\ker\Omega_B)/G \to N/G$ is a
symplectic Lie algebroid. In addition, if we have a Hamiltonian
function $H_M:M\to \R$ we obtain, under some mild hypotheses,
reduced Hamiltonian dynamics.

In the particular case when the Lie algebroid is the tangent
bundle of a symplectic manifold, our reduction procedure is just
the well known Cartan symplectic reduction in the presence of a
symmetry Lie group. This example is shown in Section 4, along with
other dif\/ferent interesting examples. A particular application of
our results is a ``symplectic description" of the Hamiltonian
reduction process by stages in the Poisson setting (see Section
4.3) and the reduction of the Lagrange top (see Section 4.4).

In the last part of the paper, we include an Appendix where we
describe how to induce, from a Lie algebroid structure on a vector
bundle $\tau _A:A\to M$ and a linear epimorphism $\pi _A:A\to
\widetilde{A}$ over $\pi _M:M\to \widetilde{M}$, a Lie algebroid
structure on $\widetilde{A}$. An equivalent dual version of this
result was proved in \cite{CNS}.

\section{Hamiltonian Mechanics and symplectic Lie algebroids}\label{section1}

\subsection{Lie algebroids}\label{section1.1}
Let $A$ be a vector bundle of $rank$ $n$ over a manifold $M$ of
dimension $m$  and $\tau_A:A\to M$ be the vector bundle
projection. Denote by $\Gamma(A)$ the $C^\infty(M)$-module of
sections of $\tau_A:A\to M$. A~{\it Lie algebroid structure }
$(\lcf\cdot,\cdot\rcf_A,\rho_A)$ on $A$ is a Lie bracket
$\lcf\cdot,\cdot\rcf_A$ on the space $\Gamma(A)$ and a bundle map
$\rho_A:A\to TM$, called {\it the anchor map}, such that if we
also denote by $\rho_A:\Gamma(A)\to {\mathfrak X}(M)$ the
homomorphism of $C^\infty(M)$-modules induced by the anchor map
then
\[
\lcf X,fY\rcf_A=f\lcf X,Y\rcf_A + \rho_A(X)(f)Y,
\]
for $X,Y\in \Gamma(A)$ and $f\in C^\infty(M)$. The triple
$(A,\lcf\cdot,\cdot\rcf_A,\rho_A)$ is called {\it a Lie algebroid
over} $M$ (see \cite{Mac}).

If $(A,\lcf\cdot,\cdot\rcf_A,\rho_A)$ is a Lie algebroid over $M,$
then the anchor map $\rho_A:\Gamma(A)\to {\mathfrak X}(M)$ is
a~homomorphism between the Lie algebras
$(\Gamma(A),\lcf\cdot,\cdot\rcf_A)$ and $({\mathfrak
X}(M),[\cdot,\cdot])$.

Trivial examples of Lie algebroids are real Lie algebras of f\/inite
dimension and the tangent bundle $TM$ of an arbitrary manifold
$M$. Other examples of Lie algebroids are the following ones:
\begin{itemize}\itemsep=0pt
\item {\bf The Lie algebroid associated with an inf\/initesimal
action.}
\end{itemize}

Let ${\mathfrak g}$ be a real Lie algebra of f\/inite dimension and
$\Phi:{\mathfrak g}\to {\mathfrak X}(M)$ an inf\/initesimal left
action of ${\mathfrak g}$ on a manifold $M$, that is, $\Phi$ is a
$\R$-linear map and
\[
\Phi([\xi,\eta]_{\mathfrak g})=-[\Phi(\xi),\Phi(\eta)], \qquad \mbox{for
all} \ \ \xi,\eta\in {\mathfrak g},
\]
where $[\cdot,\cdot]_{\mathfrak g}$ is the Lie bracket on
${\mathfrak g}.$ Then, the trivial vector bundle $\tau_A:A=M\times
{\mathfrak g}\to M$ admits a Lie algebroid structure. The anchor
map $\rho_A:A\to TM$ of $A$ is given by
\[
\rho_A(x,\xi)=-\Phi(\xi)(x), \qquad \mbox{for} \ \ (x,\xi)\in
M\times {\mathfrak g}=A.
\]

On the other hand, if $\xi$ and $\eta$ are elements of ${\mathfrak
g}$
 then $\xi$ and $\eta$ induce constant sections of $A$ which we
 will also denote by $\xi$ and $\eta$. Moreover, the Lie bracket
 $\lcf\xi,\eta\rcf_A$ of $\xi$ and $\eta$ in $A$ is the constant
 section on $A$ induced by $[\xi,\eta]_{\mathfrak g}$. In other words,
 \[
 \lcf\xi,\eta\rcf_A=[\xi,\eta]_{\mathfrak g}.
 \]
The resultant Lie algebroid is called the \emph{Lie algebroid
associated with the infinitesimal action} $\Phi$.

Let $\{\xi_\alpha\}$ be a basis of ${\mathfrak g}$ and $(x^i)$ be a
system of local coordinates on an open subset $U$ of $M$ such that
\[
[\xi_\alpha,\xi_\beta]_{\mathfrak
g}=c_{\alpha\beta}^\gamma\xi_\gamma,\qquad
\Phi(\xi_\alpha)=\Phi_\alpha^i\frac{\partial }{\partial x^i}.
\]
If $\widetilde{\xi}_\alpha:U\to M\times {\mathfrak g}$ is the map
def\/ined by
\[
\widetilde{\xi}_\alpha(x)=(x,\xi_\alpha),\qquad \mbox{for all} \ \ x\in U,
\]
then $\{\widetilde{\xi}_\alpha\}$ is a local basis of sections of
the action Lie algebroid. In addition, the corresponding local
structure functions with respect to $(x^i)$ and $\{
\tilde{\xi}_\alpha \}$ are
\[
C_{\alpha\beta}^\gamma=c_{\alpha\beta}^\gamma,\qquad
\rho_\alpha^i=\Phi_\alpha^i.
\]

\begin{itemize}
\item
 {\bf The Atiyah (gauge) algebroid associated with a principal
 bundle.}
\end{itemize}
Let $\pi_P:P\to M$ be a principal left $G$-bundle. Denote by
$\Psi:G\times P\to P$ the free action of $G$ on $P$ and by
$T\Psi:G\times TP\to TP$ the tangent action of $G$ on $TP.$ The
space $TP/G$ of orbits of the action is a vector bundle over the
manifold $M$ with vector bundle projection $\tau_{TP}|G:TP/G\to
P/G\cong M$ given by
\[
(\tau_{TP}|G)([v_p])=[\tau_{TP}(v_p)]=[p],\qquad \mbox{for}\ \ v_p\in T_pP,
\]
$\tau_{TP}:TP\to P$ being the canonical projection. A section of
the vector bundle $\tau_{TP}|G:TP/G\to P/G\cong M$ may be
identif\/ied with a vector f\/ield on $P$ which is $G$-invariant.
Thus, using that every $G$-invariant vector f\/ield on $P$ is
$\pi_P$-projectable on a vector f\/ield on $M$ and that the standard
Lie bracket of two $G$-invariant vector f\/ields is also a
$G$-invariant vector f\/ield, we may induce a Lie algebroid
structure $(\lcf\cdot,\cdot\rcf_{TP/G},\rho_{TP/G})$ on the vector
bundle $\tau_{TP}|G:TP/G\to P/G\cong M$. The Lie algebroid
$(TP/G,\lcf\cdot,\cdot\rcf_{TP/G},\rho_{TP/G})$ is called the
\emph{Atiyah (gauge) algebroid associated with the principal
bundle $\pi_P:P\to M$} (see \cite{LMM,Mac}).

Let $D:TP\to {\mathfrak g}$ be a connection in the principal
bundle $\pi_P:P\to M$ and $R:TP\oplus TP\to {\mathfrak g}$ be the
curvature of $D$. We choose a local trivialization of the
principal bundle $\pi_P:P\to M$ to be $U\times G$, where $U$ is an
open subset of $M$. Suppose that $e$ is the identity element of
$G$, that $(x^i)$ are local coordinates on $U$ and that
$\{\xi_a\}$ is a basis of ${\mathfrak g}.$

Denote by $\{\xi_a^L\}$ the corresponding left-invariant vector
f\/ields on $G$. If
\[
D\left(\frac{\partial }{\partial
x^i}_{|_{(x,e)}}\right)=D_i^a(x)\xi_a,\qquad
R\left(\frac{\partial }{\partial
x^i}_{|_{(x,e)}},\frac{\partial }{\partial x^j}_{|_{(x,e)}}\right)=
R_{ij}^a(x)\xi_a,
\]
for $x\in U$, then the horizontal lift of the vector f\/ield
$\frac{\partial }{\partial x^i}$ is the vector f\/ield
on $U\times G$ given by
\[
\left(\frac{\partial }{\partial x^i}\right)^h=\frac{\partial }{\partial
x^i}-D_i^a\xi_a^L.
\]
Therefore, the vector f\/ields on $U\times G$
$\{e_i=\frac{\partial }{\partial
x^i}-D_i^a\xi_a^L,e_b=\xi_b^L\}$ are $G$-invariant and they def\/ine
a local basis $\{e_i',e_b'\}$ of $\Gamma(TP/G).$ The corresponding
local structure functions of $\tau_{TP/G}:TP/G\to M$ are
\begin{gather}
C_{ij}^k= C_{ia}^j=-C_{ai}^j=
C_{ab}^i=0,\qquad C_{ij}^a=-R_{ij}^a,\qquad
C_{ia}^c=-C_{ai}^c=c_{ab}^cD_i^b,\qquad
C_{ab}^c=c_{ab}^c,\nonumber\\
\rho_i^j=\delta_{ij},\qquad \rho_i^a=\rho_a^i=\rho_a^b=0\label{LSF2}
\end{gather}
(for more details, see \cite{LMM}).

An important operator associated with a Lie algebroid
$(A,\lcf\cdot,\cdot\rcf_A,\rho_A)$ over a manifold $M$ is the
\emph{differential $d^A:\Gamma(\wedge^k A^*)\to
\Gamma(\wedge^{k+1}A^*)$ of $A$ } which is def\/ined as follows
\begin{gather*}
d^A \mu(X_0,\dots,
X_k)=\sum_{i=0}^k(-1)^i\rho_A(X_i)(\mu(X_0,\dots,
\widehat{X_i},\dots, X_k))\\
\phantom{d^A \mu(X_0,\dots,X_k)=}{}+ \sum_{i<
j}(-1)^{i+j}\mu(\lcf X_i,X_j\rcf_A,X_0,\dots,
\widehat{X_i},\dots,\widehat{X_j},\dots ,X_k),
\end{gather*}
for $\mu \in \Gamma(\wedge^k A^*)$ and $X_0,\dots ,X_k\in
\Gamma(A).$ It follows that $(d^A)^2=0.$ Note that if $A=TM$ then
$d^{TM}$ is the standard dif\/ferential exterior for the manifold
$M.$

On the other hand, if $(A,\lcf\cdot,\cdot\rcf_A,\rho_A)$ and
$(A',\lcf\cdot,\cdot\rcf_{A'},\rho_{A'})$ are Lie algebroids over
$M$ and $M'$, respectively, then a morphism of vector bundles
$F:A\to A'$ is a \emph{Lie algebroid morphism} if
\begin{equation}\label{Morph}
d^A(F^*\phi')= F^*(d^{A'}\phi'), \qquad \mbox{for}\ \ \phi'\in
\Gamma(\wedge^k(A')^*).
\end{equation}
Note that $F^*\phi'$ is the section of the vector bundle
$\wedge^kA^*\to M$ def\/ined by
\[
(F^*\phi')_x(a_1,\dots ,a_k)=\phi'_{f(x)}(F(a_1),\dots ,F(a_k)),
\]
for $x\in M$ and $a_1,\dots ,a_k\in A$, where $f:M\to M'$ is the
mapping associated with $F$ between~$M$ and $M'$. We remark that
(\ref{Morph}) holds if and only if
\begin{gather}
d^A(g'\circ f)=F^{*}(d^{A'}g'), \qquad \mbox{for}\ \ g'\in
C^\infty(M'),\nonumber\\
d^A(F^*\alpha')=F^*(d^{A'}\alpha'),\qquad \mbox{for} \ \ \alpha'\in
\Gamma((A')^*).\label{Morph1}
\end{gather}
 This def\/inition of a Lie algebroid
morphism is equivalent of the original one given in \cite{HM}.

If $M=M'$ and $f=id:M\to M$ then it is easy to prove that $F$ is a
Lie algebroid morphism if and only  if
\[
F\lcf X,Y\rcf_{A}=\lcf FX,FY \rcf_{A'},\qquad
\rho_{A'}(FX)=\rho_A(X),
\]
for $X,Y\in \Gamma(A).$

If $F$ is a Lie algebroid morphism, $f$ is an injective immersion
and $F$ is also injective, then the Lie algebroid
$(A,\lcf\cdot,\cdot\rcf_A,\rho_A)$ is a \emph{Lie subalgebroid} of
$(A',\lcf\cdot,\cdot\rcf_{A'},\rho_{A'})$.

\subsection{Symplectic Lie algebroids}\label{section1.2}

Let $(A,\lcf\cdot,\cdot\rcf_A,\rho_A)$ be a Lie algebroid over a
manifold $M$. Then $A$ is said to be a \emph{symplectic Lie
algebroid} if there exists a section $\Omega_A$ of the vector
bundle $\wedge^2A^*\to M$ such that $\Omega_A$ is nondegenerate
and $d^A\Omega_A=0$ (see \cite{LMM}). The f\/irst example of a
symplectic Lie algebroid is the tangent bundle of a symplectic
manifold.

It is clear that the rank of a symplectic Lie algebroid $A$ is
even. Moreover, if $f\in C^\infty(M)$ one may introduce the
\emph{Hamiltonian section ${\mathcal H}_f^{\Omega_A}\in \Gamma(A)$
of $f$ with respect to $\Omega_A$} which is characterized by the
following condition
\begin{equation}\label{HS}
i({\mathcal H}_f^{\Omega_A})\Omega _A=d^Af.
\end{equation}

On the other hand, the map $\flat_{\Omega_A}:\Gamma(A)\to
\Gamma(A^*)$ given by
\[
\flat_{\Omega_A}(X)=i(X)\Omega_A,\qquad \mbox{for}\ \ X\in
\Gamma(A),
\]
is an isomorphism of $C^\infty(M)$-modules. Thus, one may def\/ine a
section $\Pi_{\Omega_A}$ of the vector bundle $\wedge^2A\to A$ as
follows
\[
\Pi_{\Omega_A}(\alpha,\beta)=\Omega_A(\flat_{\Omega_A}^{-1}(\alpha),
\flat_{\Omega_A}^{-1}(\beta)),\qquad \mbox{for} \ \ \alpha,\beta\in
\Gamma(A^*).
\]

$\Pi_{\Omega_A}$ is a \emph{triangular matrix} for the Lie
algebroid $A$ ($\Pi_{\Omega_A}$ is an $A$-Poisson bivector f\/ield
on the Lie algebroid $A$ in the terminology of \cite{CW}) and the
pair $(A,A^*)$ is a \emph{triangular Lie bialgebroid} in the sense
of Mackenzie and Xu \cite{MX}. Therefore, the base space $M$
admits a Poisson structure, that is, a bracket of functions
\[
\{\cdot,\cdot\}_M:C^\infty(M)\times C^\infty(M)\to C^\infty(M)
\]
which satisf\/ies the following properties:
\begin{enumerate} \itemsep=0pt
\item
$\{\cdot,\cdot\}_M$ is $\R$-bilinear and skew-symmetric;
\item It
is a derivation in each argument with respect to the standard
product of functions, i.e.,
\[
\{ff',g\}_M=f\{f',g\}_M + f'\{f,g\}_M;
\]
\item It satisf\/ies the Jacobi identity, that is,
\[
\{f,\{g,h\}_M\}_M+\{g,\{h,f\}_M\}_M + \{h,\{f,g\}_M\}_M=0.
\]
\end{enumerate}
In fact, we have that
\begin{equation}\label{PB}
\{f,g\}_M=\Omega_A({\mathcal H}_f^{\Omega_A},{\mathcal
H}_g^{\Omega_A}),\qquad \mbox{for} \quad f,g\in C^\infty(M).
\end{equation}

Using the Poisson bracket $\{\cdot,\cdot\}_M$ one may consider the
\emph{Hamiltonian vector field} of a real function $f\in
C^\infty(M)$ as the vector f\/ield ${\mathcal
H}_f^{\{\cdot,\cdot\}_M}$ on $M$ def\/ined by
\[
{\mathcal H}_f^{\{\cdot,\cdot\}_M}(g)=-\{f,g\}_M,\,\qquad \mbox{for all} \ \ g\in C^\infty(M).\]

From (\ref{HS}) and (\ref{PB}), we deduce that $\rho_A({\mathcal
H}_f^{\Omega_A})={\mathcal H}_f^{\{\cdot,\cdot\}_M}$. Moreover,
the f\/low of the vector f\/ield ${\mathcal H}_f^{\{\cdot,\cdot\}_M}$
on $M$ is covered by a group of $1$-parameter automorphisms of
$A$: it is just the (inf\/initesimal f\/low) of the $A$-vector f\/ield
${\mathcal H}_f^{\Omega_A}$ (see~\cite{R}).

 A symplectic Lie algebroid may be associated
with an arbitrary Lie algebroid as follows (see~\cite{LMM}).

Let $(A,\lcf\cdot,\cdot\rcf_A,\rho_A)$ be a Lie algebroid of rank
$n$ over a manifold $M$ and $\tau_{A^*}:A^*\to M$ be the vector
bundle projection of the dual bundle $A^*$ to $A$. Then, we
consider the $A$-tangent bundle to $A^*$ as the subset of $A\times
TA^*$ given by
\[
{\mathcal T}^AA^*=\{(b,v)\in A\times
TA^*/\rho_A(b)=(T\tau_{A^*})(v)\}.
\]
${\mathcal T}^AA^*$ is a vector bundle over $A^*$ of rank $2n$ and
the vector bundle projection $\tau_{{\mathcal T}^AA^*}:{\mathcal
T}^AA^*\to A^*$ is def\/ined by
\[
\tau_{{\mathcal T}^AA^*}(b,v)=\tau_{TA^*}(v),\qquad \mbox{for} \ \ (b,v)\in
{\mathcal T}^AA^*.
\]

A section $\widetilde{X}$ of $\tau_{{\mathcal T}^AA^*}: {\mathcal
T}^AA^* \to A^*$ is said to be \emph{projectable} if there exists
a section $X$ of $\tau_A: A \to M$ and a vector f\/ield $S \in
{\mathfrak X}(A^*)$ which is $\tau_{A^*}$-projectable to the
vector f\/ield $\rho_A(X)$ and such that $\tilde{X}(p) =
(X(\tau_{A^*}(p)), S(p))$, for all $p \in A^*$. For such a
projectable section $\tilde{X}$, we will use the following
notation $\tilde{X} \equiv (X, S)$. It is easy to prove that one
may choose a local basis of projectable sections of the space
$\Gamma({\mathcal T}^AA^*)$.

The vector bundle $\tau_{{\mathcal T}^AA^*}: {\mathcal T}^AA^* \to
A^*$ admits a Lie algebroid structure $(\lcf \cdot , \cdot
\rcf_{{\mathcal T}^AA^*} , \rho_{{\mathcal T}^AA^* })$. Indeed, if
$(X, S)$ and $(Y, T)$ are projectable sections then
\[
\lcf (X, S), (Y, T) \rcf_{{\mathcal T}^AA^* }= (\lcf X, Y\rcf _A,
[S, T]), \qquad  \rho_{{\mathcal T}^AA^*} (X, S) = S.
\]

$({\mathcal T}^AA^*,\lcf\cdot,\cdot\rcf_{{\mathcal T}^AA^*}
,\rho_{{\mathcal T}^AA^*} )$ is the \emph{$A$-tangent bundle to
$A^*$} or \emph{the prolongation of $A$ over the fibration
$\tau_{A^*}:A^*\to M$} (for more details, see \cite{LMM}).

Moreover, one may introduce a canonical section $\Theta_{{\mathcal
T}^AA^*}$ of the vector bundle $({\mathcal T}^AA^*)^*\to A^*$ as
follows
\[
\Theta_{{\mathcal T}^AA^*}(\gamma)(b,v)=\gamma(b),
\]
for $\gamma\in A^*$ and $(b,v)\in {\mathcal T}_\gamma^AA^*$.
$\Theta_{{\mathcal T}^AA^*}$ is called the \emph{Liouville section
associated with} $A$ and $\Omega_{{\mathcal T}^AA^*}=-d^{{\mathcal
T}^AA^*}\Theta_{{\mathcal T}^AA^*}$ is the \emph{canonical
symplectic section associated with $A$}. $\Omega_{{\mathcal
T}^AA^*}$ is a~symplectic section for the Lie algebroid ${\mathcal
T}^AA^*$.

Therefore, the base space $A^*$ admits a Poisson structure
$\{\cdot,\cdot\}_{A^*}$ which is characterized by the following
conditions
\begin{gather*}
\{f\circ \tau_{A^*},g\circ \tau_{A^*}\}_{A^*}=0,\qquad \{f\circ
\tau_{A^*},\widehat{X}\}_{A^*}=\rho_A(X)(f)\circ \tau_{A^*},\\
\{\widehat{X},\widehat{Y}\}_{A^*}=-\widehat{\lcf X,Y\rcf_A},
\end{gather*}
for $f,g\in C^\infty(M)$ and $X,Y\in \Gamma(A).$ Here,
$\widehat{Z}$ denotes the linear function on $A^*$ induced by a~section $Z\in \Gamma(A).$ The Poisson structure
$\{\cdot,\cdot\}_{A^*}$ is called the \emph{canonical linear
Poisson structure on $A^*$ associated with the Lie algebroid $A$.}

Next, suppose that $(x^i)$ are local coordinates on an open subset
$U$ of $M$ and that $\{e_\alpha\}$ is a~local basis of $\Gamma(A)$
on $U$. Denote by $(x^i,y_\alpha)$ the corresponding local
coordinates on $A^*$ and by $\rho_\alpha^i$,
$C_{\alpha\beta}^\gamma$ the local structure functions of $A$ with
respect to the coordinates $(x^i)$ and to the basis
$\{e_\alpha\}$. Then, we may consider the local sections
$\{{\mathcal X}_\alpha,{\mathcal P}^\alpha\}$ of the vector bundle
$\tau_{{\mathcal T}^AA^*}:{\mathcal T}^AA^*\to A^*$ given by
\[
{\mathcal X}_\alpha=\left(e_\alpha\circ
\tau_{A^*},\rho_\alpha^i\frac{\partial }{\partial x^i}\right),\qquad
{\mathcal P}^\alpha=\left(0,\frac{\partial }{\partial y_\alpha}\right).
\]
We have that $\{{\mathcal X}_\alpha,{\mathcal P}^\alpha\}$ is a
local basis of $\Gamma({\mathcal T}^AA^*)$ and
\begin{gather*}
\lcf {\mathcal X}_\alpha,{\mathcal X}_\beta\rcf_{{\mathcal
T}^AA^*}=C_{\alpha\beta}^\gamma{\mathcal X}_\gamma,\qquad
\lcf{\mathcal X}_\alpha,{\mathcal P}^\beta\rcf_{{\mathcal
T}^AA^*}=\lcf{\mathcal P}^\alpha,{\mathcal P}^\beta\rcf_{{\mathcal
T}^AA^*}=0,
\\
\rho_{{\mathcal T}^AA^*}({\mathcal
X}_\alpha)=\rho_\alpha^i\frac{\partial }{\partial x^i},\qquad
\rho_{{\mathcal T}^AA^*}({\mathcal P}^\alpha)=\frac{\partial
}{\partial y_\alpha},
\\
\Theta_{{\mathcal T}^AA^*}=y_\alpha{\mathcal X}^\alpha,\qquad
\Omega_{{\mathcal T}^AA^*}={\mathcal X}^\alpha\wedge {\mathcal
P}_\alpha + \frac{1}{2}C_{\alpha\beta}^\gamma y_\gamma{\mathcal
X}^\alpha\wedge {\mathcal X}^\beta,
\end{gather*}
$\{{\mathcal X}^\alpha,{\mathcal P}_\alpha\}$ being the dual basis
of $\{{\mathcal X}_\alpha,{\mathcal P}^\alpha\}$. Moreover, if
$H,H':A^*\to \R$ are real functions on~$A^*$ it follows that
\begin{gather} \{H,H'\}_{A^*}=-\frac{\partial H}{\partial
y_\alpha}\frac{\partial H'}{\partial
y_\beta}C_{\alpha\beta}^\gamma y_\gamma + \left(
\frac{\partial H}{\partial
x^i}\frac{\partial H'}{\partial y_\alpha}
-\frac{\partial H}{\partial
y_\alpha}\frac{\partial  H'}{\partial
x^i}\right)\rho_\alpha^i,\nonumber\\
{\mathcal H}_H^{\Omega_{{\mathcal
T}^AA^*}}=\frac{\partial H}{\partial
y_\alpha}{\mathcal X}_\alpha-\left(
\rho_\alpha^i \frac{\partial H}{\partial x^i}
+C_{\alpha\beta}^\gamma y_\gamma  \frac{\partial
H}{\partial y_\beta}\right ) {\mathcal P}^\alpha,\label{Hamilton}\\
{\mathcal
H}_H^{\{\cdot,\cdot\}_{A^*}}= \frac{\partial
H}{\partial y_\alpha}\rho_\alpha^i \frac{\partial
}{\partial x^i} - \left( \rho_\alpha^i \frac{\partial
H}{\partial x^i} + C_{\alpha\beta}^\gamma y_\gamma
\displaystyle\frac{\partial H}{\partial y_\beta}\right)
\frac{\partial }{\partial y_\alpha},\nonumber
\end{gather}
(for more details, see \cite{LMM}).

\begin{example} \qquad{}\null
\begin{enumerate}\itemsep=0pt
\item If $A$ is the standard Lie algebroid $TM$ then ${\mathcal
T}^AA^*=T(T^*M),$ $\Omega_{{\mathcal T}^AA^*}$ is the canonical
symplectic structure on $A^*=T^*M$ and $\{\cdot,\cdot\}_{T^*M}$ is
the canonical Poisson bracket on $T^*M$ induced by
$\Omega_{{\mathcal T}^AA^*}$.

\item If $A$ is the Lie algebroid associated with an inf\/initesimal
action $\Phi:{\mathfrak g}\to {\mathfrak X}(M)$ of a Lie algebra
${\mathfrak g}$ on a manifold $M$ (see Section \ref{section1.1}),
then the Lie algebroid ${\mathcal T}^AA^*\to A^*$ may be
identif\/ied with the trivial vector bundle
\[
(M\times {\mathfrak g}^*)\times ({\mathfrak g}\times {\mathfrak
g}^*)\to M\times {\mathfrak g}^*
\]
and, under this identif\/ication, the canonical symplectic section
$\Omega_{{\mathcal T}^AA^*}$ is given by
\[
\Omega_{{\mathcal
T}^AA^*}(x,\alpha)((\xi,\beta),(\xi',\beta'))=\beta'(\xi)-\beta(\xi')
+ \alpha[\xi,\xi']_{\mathfrak g},
\]
for $(x,\alpha)\in A^*=M\times {\mathfrak g}^*$ and
$(\xi,\beta),(\xi',\beta')\in {\mathfrak g}\times {\mathfrak
g}^*$, where $[\cdot,\cdot]_{\mathfrak g}$ is the Lie bracket on
${\mathfrak g}.$

The anchor map $\rho _{{\mathcal T}^AA^*}:(M\times \mathfrak
g^*)\times (\mathfrak g\times \mathfrak g^*)\to T(M\times
\mathfrak g^*)\cong TM\times (\mathfrak g^*\times \mathfrak g^*)$
of the $A$-tangent bundle to $A^*$ is
\[
\rho _{{\mathcal T}^AA^*}((x,\alpha),(\xi,\beta))=(-\Phi (\xi
)(x),\alpha, \beta)
\]
and the Lie bracket of two constant sections $(\xi ,\beta)$ and
$(\xi ',\beta')$ is just the constant section $([\xi ,\xi
']_\mathfrak g,0)$.

Note that in the particular case when $M$ is a single point  (that
is, $A$ is the real algebra~${\mathfrak g}$) then the linear
Poisson bracket $\{\cdot,\cdot\}_{{\mathfrak g}^*}$ on ${\mathfrak
g}^*$ is just the (minus) Lie--Poisson bracket induced by the Lie
algebra ${\mathfrak g}.$

\item Let $\pi_P:P\to M$ be a principal $G$-bundle and $A=TP/G$ be
the Atiyah algebroid associated with $\pi_P:P\to M$. Then, the
cotangent bundle $T^*P$ to $P$ is the total space of a principal
$G$-bundle over $A^*\cong T^*P/G$ and the Lie algebroid ${\mathcal
T}^AA^*$ may be identif\/ied with the Atiyah algebroid $T(T^*P)/G\to
T^*P/G$ associated with this principal $G$-bundle (see
\cite{LMM}). Moreover, the canonical symplectic structure on
$T^*P$ is $G$-invariant and it induces a symplectic section
$\widetilde{\Omega}$ on the Atiyah algebroid $T(T^*P)/G\to
T^*P/G.$ $\widetilde{\Omega}$ is just the canonical symplectic
section of ${\mathcal T}^AA^*\cong T(T^*P)/G\to A^*\cong T^*P/G.$
Finally, the linear Poisson bracket on $A^*\cong T^*P/G$ is
characterized by the following property:  if on $T^*P$ we consider
the Poisson bracket induced by the canonical symplectic structure
then the canonical projection $\pi_{T^*P}:T^*P\to T^*P/G$ is a
Poisson morphism.

We remark that, using a connection in the principal $G$-bundle
$\pi_P:P\to M$, the space $A^*\cong T^*P/G$ may be identif\/ied with
the Whitney sum $W=T^*M\oplus_M\widetilde{\mathfrak g}^*,$ where
$\widetilde{\mathfrak g}^*$ is the coadjoint bundle associated
with the principal bundle $\pi_P:P\to M.$ In addition, under the
above identif\/ication, the Poisson bracket on $A^*\cong T^*P/G$ is
just the so-called \emph{Weinstein space Poisson bracket} (for
more details, see \cite{OR}).
\end{enumerate}
\end{example}

\subsection{Hamilton equations and symplectic Lie algebroids}\label{section1.3}
Let $A$ be a Lie algebroid over a manifold $M$ and $\Omega_A$ be a
symplectic section of $A$. Then, as we know, $\Omega_A$ induces a
Poisson bracket $\{\cdot,\cdot\}_M$ on $M$.

A Hamiltonian function for $A$ is a real $C^\infty$-function
$H:M\to \R$ on $M$. If $H$ is a Hamiltonian function one may
consider the Hamiltonian section ${\mathcal H}^{\Omega_A}_H$ of
$H$ with respect to $\Omega_A$ and the Hamiltonian vector f\/ield
${\mathcal H}_H^{\{\cdot,\cdot\}_M}$ of $H$ with respect to the
Poisson bracket $\{\cdot,\cdot\}_M$ on $M$.

\emph{The solutions of the Hamilton equations for $H$ in $A$} are
just the integral curves of the vector f\/ield ${\mathcal
H}_H^{\{\cdot,\cdot\}_M}$.

Now, suppose that our symplectic Lie algebroid is the $A$-tangent
bundle to $A^*$, ${\mathcal T}^AA^*,$ where $A$ is an arbitrary
Lie algebroid. As we know, the base space of ${\mathcal T}^AA^*$
is $A^*$ and the corresponding Poisson bracket
$\{\cdot,\cdot\}_{A^*}$ on $A^*$ is just the canonical linear
Poisson bracket on $A^*$ associated with the Lie algebroid $A$.
Thus, if $H:A^*\to \R$ is a Hamiltonian function for ${\mathcal
T}^AA^*,$ we have that the solutions of the Hamilton equations for
$H$ in ${\mathcal T}^AA^*$ (or simply, the solutions of the
Hamilton equations for $H$) are the integral curves of the
Hamiltonian vector f\/ield ${\mathcal H}_H^{\{\cdot,\cdot\}_{A^*}}$
(see \cite{LMM}).

If $(x^i)$ is a system of local coordinates on an open subset $U$
of $M$ and $\{e_\alpha\}$ is a local basis of~$\Gamma(A)$ on $U$,
we may consider the corresponding system  of local coordinates
$(x^i,y_\alpha)$ on $\tau_{A^*}^{-1}(U)\subseteq A^*$. In
addition, using (\ref{Hamilton}), we deduce that a curve $t\to
(x^i(t),y_\alpha(t))$ on $\tau_{A^*}^{-1}(U)$ is a~solution of the
Hamilton equations for $H$ if and only if
\begin{equation}\label{HamEq}
\frac{dx^i}{dt}=\rho_\alpha^i\frac{\partial H}{\partial
x^i},\qquad \frac{dy_\alpha}{dt}=-\left(
C_{\alpha\beta}^\gamma y_\gamma\frac{\partial H}{\partial y_\beta} +
\rho_\alpha^i\frac{\partial H}{\partial x^i} \right),
\end{equation}
(see \cite{LMM}).

\begin{example}\qquad {}
\begin{enumerate}\itemsep=0pt
\item Let $A$ be the Lie algebroid associated with an
inf\/initesimal action $\Phi:{\mathfrak g}\to {\mathfrak X}(M)$ of a
Lie algebra ${\mathfrak g}$ on the manifold $M$. If $H:A^*=M\times
{\mathfrak g}^*\to \R$ is a Hamiltonian function, $\{\xi_\alpha\}$
is a basis of ${\mathfrak g }$ and $(x^i)$ is a system of local
coordinates on $M$ such that
\[
[\xi_\alpha,\xi_\beta]_{\mathfrak
g}=c_{\alpha\beta}^\gamma\xi_\gamma,\qquad
\Phi(\xi_\alpha)=\Phi_\alpha^i\frac{\partial }{\partial x^i},
\]
then the curve $t\to (x^i(t),y_\alpha(t))$ on $A^*=M\times
{\mathfrak g^*}$ is a solution of the Hamilton equations for $H$
if and only if
\[
\frac{dx^i}{dt}=\Phi_\alpha^i\frac{\partial H}{\partial
x^i},\qquad \frac{dy_\alpha}{dt}=-\left(
c_{\alpha\beta}^\gamma y_\gamma\frac{\partial H}{\partial y_\beta} +
\Phi_\alpha^i\frac{\partial H}{\partial x^i}\right).
\]

Note that if $M$ is a single point (that is, $A$ is the Lie
algebra ${\mathfrak g}$) then the above equations are just
\emph{the (minus) Lie--Poisson equations on ${\mathfrak g}^*$ for
the Hamiltonian function $H:{\mathfrak g}^*\to \R.$}

\item Let $\pi_P:P\to M$ be a principal $G$-bundle, $D:TP\to
{\mathfrak g}$ be a principal connection and $R:TP\oplus TP\to
{\mathfrak g}$ be the curvature of $D$. We choose a local
trivialization of the principal bundle $\pi_P:P\to M$ to be
$U\times G$, where $U$ is an open subset of $M$. Suppose that
$(x^i)$ are local coordinates on $U$ and that $\{\xi_a\}$ is a
basis of ${\mathfrak g}$ such that $[\xi _a,\xi _b]_\mathfrak
g=c^c_{ab}\xi _c$. Denote by $D_i^a$ (respectively, $R_{ij}^a$)
the components of $D$ (respectively, $R$) with respect to the
coordinates~$(x^i)$ and to the basis $\{\xi_a\}$ and by
$(x^i,y_\alpha=p_i,\bar{p}_a)$ the corresponding local coordinates
on $A^*\cong T^*P/G$ (see Section \ref{section1.1}). If
$h:T^*P/G\to \R$ is a Hamiltonian function and $c:t\to
(x^i(t),p_i(t),\bar{p}_a(t))$ is a curve on $A^*\cong T^*P/G$
then, using (\ref{LSF2}) and (\ref{HamEq}), we conclude that $c$
is a solution of the Hamilton equations for $h$ if and only if
\begin{gather}
 \frac{dx^i}{dt}= \frac{\partial
h}{\partial
p_i},\qquad \frac{dp_i}{dt}=- \frac{\partial
h}{\partial x^i} + R_{ij}^a\bar{p}_a \frac{\partial
h}{\partial p_j}-c_{ab}^cD_i^b\bar{p}_c \frac{\partial
h}{\partial
\bar{p}_a},\nonumber\\
\frac{d\bar{p}_a}{dt} =c_{ab}^cD_i^b\bar{p}_c \frac{\partial
h}{\partial p_i}-c_{ab}^c\bar{p}_c \frac{\partial
h}{\partial \bar{p}_b}.\label{HamPo}
\end{gather}
These equations are just \emph{the Hamilton--Poincar{\'e} equations
associated with the $G$-invariant Hamiltonian function $H=h\circ
\pi_{T^*P},$} $\pi_{T^*P}:T^*P\to T^*P/G$ being the canonical
projection (see \cite{LMM}).

Note that in the particular case when $G$ is the trivial Lie group,
equations (\ref{HamPo}) reduce to
\[
\frac{dx^i}{dt}=\frac{\partial h}{\partial p_i},\qquad
\frac{dp_i}{dt}=-\frac{\partial h}{\partial x^i},
\]
which are the classical Hamilton equations for the Hamiltonian
function \mbox{$h\!:\!T^*P{\cong} T^*M\!\to\! \R$}.

\end{enumerate}

\end{example}

\section{Reduction of the Hamiltonian  dynamics
on a symplectic Lie\\ algebroid by a Lie subalgebroid and a symmetry Lie group}\label{section2}

\subsection{Reduction of the symplectic Lie
algebroid}\label{section2.1}

Let $A$ be a Lie algebroid over a manifold $M$ and $\Omega_A$ be a
symplectic section of $A.$ In addition, suppose that $B$ is a Lie
subalgebroid over the submanifold $N$ of $M$ and denote by
$i_B:B\to A$ the corresponding monomorphism of Lie algebroids. The
following commutative diagram illustrates the above situation

\begin{picture}(375,90)(40,20)
\put(195,20){\makebox(0,0){$N$}}
\put(245,25){$i_N$}\put(210,20){\vector(1,0){80}}
\put(300,20){\makebox(0,0){$M$}} \put(180,50){$\tau_B$}
\put(195,70){\vector(0,-1){40}} \put(310,50){$\tau_{A}$}
\put(300,70){\vector(0,-1){40}} \put(195,80){\makebox(0,0){$B$}}
\put(245,85){$i_B$}\put(210,80){\vector(1,0){80}}
\put(300,80){\makebox(0,0){$A$}} \end{picture}

\vspace{10pt}

$\Omega_A$ induces, in a natural way, a section $\Omega_B$ of the
real vector bundle $\wedge^2B^*\to N$. In fact,
\[
\Omega_B=i_B^*\Omega_A.
\]
Since $i_B$ is a Lie algebroid morphism, it follows that
\begin{equation}\label{closed}
d^B\Omega_B=0.
\end{equation}
However, $\Omega_B$ is not, in general, nondegenerate. In other
words, $\Omega_B$ is a presymplectic section of the Lie
subalgebroid $\tau_B:B\to N$.

For every point $x$ of $N$, we will denote by $\ker\Omega_B(x)$
the vector subspace of $B_x$ def\/ined by
\[
\ker\Omega_B(x)=\{a_x\in B_x/i(a_x)\Omega_B(x)=0\}\subseteq B_x.\]
In what follows, \emph{we will assume that the dimension of the
subspace $\ker\Omega_B(x)$ is constant, for all $x\in N.$} Thus,
the space $B/\ker\Omega_B$ is a quotient vector bundle over the
submanifold $N$.

Moreover, we may consider the vector subbundle $\ker \Omega_B$ of
$B$ whose f\/iber over the point $x\in N$ is the vector space $\ker
\Omega_B(x)$. In addition, we have that:
\begin{lemma}
$\ker \Omega_B$ is a Lie subalgebroid over $N$ of $B$.
\end{lemma}
\begin{proof} It is suf\/f\/icient to prove that
\begin{equation}\label{IntKer}
X,Y\in \Gamma(\ker\Omega_B)\Rightarrow \lcf X, Y\rcf_B\in
\Gamma(\ker \Omega_B),
\end{equation}
where $\lcf\cdot,\cdot\rcf_B$ is the Lie bracket on $\Gamma(B).$

Now, using (\ref{closed}), it follows that (\ref{IntKer}) holds.
\end{proof}

Next, we will suppose that there is \emph{a proper and free
action  $\Psi:G\times B\to B$ of a Lie group~$G$ on $B$ by vector
bundle automorphisms}. Then, the following conditions are
satisf\/ied:

\begin{enumerate}\itemsep=0pt
\item[C1)]$N$ is the total space of a principal $G$-bundle over
$\widetilde{N}$ with principal bundle projection $\pi_N:N\to
\widetilde{N}\cong N/G$ and we will denote by $\psi:G\times N\to
N$ the corresponding free action of $G$ on $N$. \item[C2)]
$\Psi_g:B\to B$ is a vector bundle isomorphism over $\psi_g:N\to
N$, for all $g\in G.$
\end{enumerate}

The following commutative diagram illustrates the above situation

\begin{picture}(375,90)(40,20)
\put(195,20){\makebox(0,0){$N$}}
\put(245,25){$\psi_g$}\put(210,20){\vector(1,0){80}}
\put(300,20){\makebox(0,0){$N$}} \put(180,50){$\tau_B$}
\put(195,70){\vector(0,-1){40}} \put(310,50){$\tau_{B}$}
\put(300,70){\vector(0,-1){40}} \put(195,80){\makebox(0,0){$B$}}
\put(245,85){$\Psi_g$}\put(210,80){\vector(1,0){80}}
\put(300,80){\makebox(0,0){$B$}}
\put(200,15){\vector(1,-1){40}}\put(295,15){\vector(-1,-1){40}}
\put(250,-30){\makebox(0,0){$\widetilde{N}\cong N/G$}}
\put(295,-10){$\pi_N$}\put(190,-10){$\pi_N$}
\end{picture}

\vspace{60pt}

In such a case, the action $\Psi$ of the Lie group $G$ on the
presymplectic Lie algebroid $(B,\Omega_B)$ is said to be \emph{presymplectic} if
\begin{equation}\label{Presyac}
\Psi_g^*\Omega_B=\Omega_B,\qquad \mbox{for all} \ \ g\in G.
\end{equation}

Now, we will prove the following result.

\begin{theorem}\label{t2.3}
Let $\Psi:G\times B\to B$ be a presymplectic action of the Lie
group $G$ on the presymplectic Lie algebroid $(B,\Omega_B)$. Then,
the action of $G$ on $B$ induces an action of $G$ on the quotient
vector bundle $B/\ker\Omega_B$ such that:
\begin{enumerate}\itemsep=0pt
\item[a)] The space of orbits $\widetilde{B}=(B/\ker\Omega_B)/G$
of this action is a vector bundle over $\widetilde{N}=N/G$ and the
diagram

\begin{picture}(375,90)(90,20)
\put(195,20){\makebox(0,0){$N$}}
\put(245,25){$\pi_N$}\put(210,20){\vector(1,0){80}}
\put(320,20){\makebox(0,0){$\widetilde{N}=N/G$}}
\put(180,50){$\tau_B$} \put(195,70){\vector(0,-1){40}}
\put(310,50){$\tau_{\widetilde{B}}$}
\put(300,70){\vector(0,-1){40}} \put(195,80){\makebox(0,0){$B$}}
\put(245,85){$\widetilde\pi_{B}$}\put(210,80){\vector(1,0){80}}
\put(340,80){\makebox(0,0){$\widetilde{B}=(B/\ker\Omega_B)/G$}}
\end{picture}

\vspace{10pt}

\noindent defines an epimorphism of vector bundles, where
$\widetilde{\pi}_{B}:B\to \widetilde{B}$ is the canonical
projection. \item[b)] There exists a unique section
$\Omega_{\widetilde{B}}$ of $\wedge^2 \widetilde{B}^*\to
\widetilde{N}$ such that
\[
 \widetilde{\pi}_B^*\Omega_{\widetilde{B}}={\Omega_B}.
 \]
\item[c)] $\Omega_{\widetilde{B}}$ is nondegenerate.
\end{enumerate}
\end{theorem}

\begin{proof} If $g\in G$ and $x\in M$ then $\Psi_g:B_x\to
 B_{\psi_g(x)}$ is a linear isomorphism and
 \[
 \Psi_g(\ker\Omega_B(x))=\ker\Omega_B(\psi_g(x)).
 \]
Thus, $\Psi$ induces an action $\widetilde{\Psi}:G\times
(B/\ker\Omega_B)\to (B/\ker\Omega_B)$ of $G$ on $B/\ker\Omega_B$
and $\widetilde{\Psi}_g$ is a vector bundle isomorphism, for all
$g\in G.$ Moreover, the canonical projection $\pi_B:B\to
B/\ker\Omega_B$ is equivariant with respect to the actions $\Psi$
and $\widetilde{\Psi}$.

On the other hand, the vector bundle projection
$\tau_{B/\ker\Omega_B}:B/\ker \Omega_B\to N$ is also equivariant
with respect to the action $\widetilde{\Psi}.$ Consequently, it
induces a smooth map
$\tau_{\widetilde{B}}:\widetilde{B}=(B/\ker\Omega_B)/G\to
\widetilde{N}=N/G$ such that the following diagram is commutative

\begin{picture}(375,90)(40,20)
\put(165,20){\makebox(0,0){$\widetilde{B}=(B/\ker\Omega_B)/G$}}
\put(250,25){$\tau_{\widetilde{B}}$}\put(215,20){\vector(1,0){75}}
\put(320,20){\makebox(0,0){$\widetilde{N}=N/G$}}
\put(150,50){$\pi_{B/\ker\Omega_B}$}
\put(195,70){\vector(0,-1){40}} \put(310,50){$\pi_N$}
\put(300,70){\vector(0,-1){40}}
\put(185,80){\makebox(0,0){$B/\ker\Omega_B$}}
\put(230,85){$\tau_{B/\ker\Omega_B}$}\put(210,80){\vector(1,0){80}}
\put(300,80){\makebox(0,0){$N$}} \end{picture}

\vspace{20pt}

Now, if $x\in N$ then, using that the action $\psi$ is free, we
deduce that the map $\pi_{B_x/\ker \Omega_B(x)}:B_x/\ker
\Omega_B(x)\to \tau_{\widetilde{B}}^{-1}(\pi_N(x))$ is bijective.
Thus, one may introduce a vector space structure on
$\tau_{\widetilde{B}}^{-1}(\pi_N(x))$ in such a way that the map
$\pi_{B_x/\ker\Omega_B(x)}:B_x/\ker\Omega_B(x)\to
\tau_{\widetilde{B}}^{-1}(\pi_N(x))$ is a linear isomorphism.
Furthermore, if $\pi_N(y)=\pi_N(x)$ then
$\pi_{B_y/\ker\Omega_B(y)}=\pi_{B_x/\ker\Omega_B(x)}\circ
 \widetilde{\Psi}_g^{-1}.$

Therefore, $\widetilde{B}$ is a vector bundle over $\widetilde{N}$
with vector bundle projection
$\tau_{\widetilde{B}}:\widetilde{B}\to \widetilde{N}$ and if $x\in
N$ then the f\/iber of $\widetilde{B}$ over $\pi_N(x)$ is isomorphic
to the vector space $B_x/\ker\Omega_B(x)$. This proves $a)$.

On the other hand, it is clear that the section $\Omega_B$ induces
a section $\Omega_{(B/\ker\Omega_B)}$ of the vector bundle
$\wedge^2(B/\ker \Omega_B)\to N$ which is characterized by the
condition
\begin{equation}\label{Omegaco}
\pi_B^*(\Omega_{(B/\ker\Omega_B)})=\Omega_B.
\end{equation}
Then, using (\ref{Presyac}) and (\ref{Omegaco}), we deduce that
the section $\Omega_{(B/\ker\Omega_B)}$ is $G$-invariant, that is,
 \[
 \widetilde{\Psi}_g^*(\Omega_{(B/\ker\Omega_B)})=\Omega_{(B/\ker\Omega_B)},\qquad
 \mbox{for all} \ \ g\in G.
 \]
Thus, there exists a unique section $\Omega_{\widetilde{B}}$ of
the vector bundle $\wedge^2\widetilde{B}^*\to \widetilde{N}$ such
that
\begin{equation}\label{Omegatil}
\pi_{B/\ker\Omega_B}^*(\Omega_{\widetilde{B}})=\Omega_{(B/\ker\Omega_B)}.
\end{equation}
This proves $b)$ (note that
$\widetilde{\pi}_B=\pi_{B/\ker\Omega_B}\circ \pi_B$).

Now, using (\ref{Omegaco}) and the fact
$\ker\pi_B=\ker\Omega_{B},$ it follows that the section
$\Omega_{(B/\ker\Omega_B)}$ is nondegenerate. Therefore, from
(\ref{Omegatil}), we conclude that $\Omega_{\widetilde{B}}$ is
also nondegenerate (note that if $x\in M$ then
$\pi_{B/\ker\Omega_B}:B_x/\ker\Omega_B(x)\to
\widetilde{B}_{\pi_N(x)}$ is a linear isomorphism).
\end{proof}

 Next, we will describe the space of sections of $\widetilde{B}=(B/\ker\Omega_B)/G$. For
 this purpose, we will use some results contained in the Appendix
 of this paper.

 Let $\Gamma(B)^p_{\widetilde{\pi}_B}$ be the space of
 $\widetilde{\pi}_B$-projectable sections of the vector
 bundle $\tau_B:B\to N.$ As we know (see the Appendix), a section
 $X$ of $\tau_B:B\to N$ is said to be
 $\widetilde{\pi}_B$-projectable if there exists a section
 $\widetilde{\pi}_B(X)$ of the vector bundle
 $\tau_{\widetilde{B}}:\widetilde{B}\to \widetilde{N}$ such that
 $\widetilde{\pi}_B\circ X=\widetilde{\pi}_B(X)\circ
 \pi_N.$ Thus, it is clear that $X$ is
 $\widetilde{\pi}_B$-projectable if and only if
 $\pi_B\circ X$ is a $\pi_{(B/\ker\Omega_B)}$-projectable
 section.

 On the other hand, it is easy to prove that the section
 $\pi_B\circ X$ is $\pi_{(B/\ker\Omega_B)}$-projectable if
 and only if it is $G$-invariant. In other words, $(\pi_B\circ X)$ is
$\pi_{(B/\ker\Omega_B)}$-projectable if and only if for every
$g\in G$ there exists $Y_g\in \Gamma(\ker\pi_B)$ such that $
\Psi_g\circ X=(X+Y_g)\circ \psi_g,$ where $\psi:G\times N\to N$ is
the corresponding free action of $G$ on $N$. Note that
 \begin{equation}\label{kernels}
 \ker\widetilde{\pi}_B=\ker\pi_B=\ker\Omega_B
 \end{equation}
 and therefore, the above facts, imply that
 \begin{equation}\label{Projectable}
 \Gamma(B)^p_{\widetilde{\pi}_B}=\{X\in \Gamma(B)/\forall\,
 g\in G,\; \exists \, Y_g\in \Gamma(\ker\Omega_B)\mbox{  and } \Psi_g\circ
 X=(X+Y_g)\circ \psi_g\}.
 \end{equation}

 Consequently, using some results of the Appendix (see (\ref{A1})
 in the Appendix), we deduce that
 \[
 \Gamma(\widetilde{B})\cong \frac{\{X\in \Gamma(B)/\forall g\in
 G,\;\; \exists Y_g\in \Gamma(\ker\Omega_B) \mbox{ and }
 \Psi_g\circ X=(X+Y_g)\circ \psi_g\}}{\Gamma(\ker\Omega_B)}
 \]
as $C^\infty(\tilde{N})$-modules.

Moreover, we may prove the following result

\begin{theorem}[{\bf The reduced symplectic Lie algebroid}]\label{t2.4} Let
$\Psi:G\times B\to B$ be a presympletic action of the Lie group
$G$ on the Lie algebroid $(B,\Omega_B)$. Then, the reduced vector
bundle $\tau_{\widetilde{B}}:\widetilde{B}=(B/\ker\Omega_B)/G\to
\widetilde{N}=N/G$ admits a unique Lie algebroid structure such
that $\widetilde{\pi}_B:B\to \widetilde{B}$ is a Lie algebroid
epimorphism if and only if the following conditions hold:
\begin{enumerate}\itemsep=0pt
\item[i)] The space $\Gamma(B)_{\widetilde{\pi}_B}^p$ is a Lie
subalgebra of the Lie algebra $(\Gamma(B),\lcf\cdot,\cdot\rcf_B)$.
\item[ii)] $\Gamma(\ker\Omega_B)$ is an ideal of this Lie
subalgebra.
\end{enumerate}

Furthermore, if the conditions $i)$ and $ii)$ hold, we get that
there exists a short exact sequence of Lie algebroids
\[
0\to\ker\Omega_B\to B\to \widetilde{B}\to 0
\]
and  that the nondegenerate section $\Omega_{\widetilde{B}}$
induces a symplectic structure on the Lie algebroid
$\tau_{\widetilde{B}}:\widetilde{B}\to \widetilde{N}.$
\end{theorem}

\begin{proof} The f\/irst part of the Theorem follows from (\ref{kernels}),
(\ref{Projectable}) and Theorem \ref{tA1} in the Appendix.

On the other hand, if the conditions $i)$ and $ii)$ in the Theorem
hold then, using Theorem \ref{t2.3} and the fact that
$\widetilde{\pi}_B$ is a Lie algebroid epimorphism, we obtain that
\[
\widetilde{\pi}_B^*(d^{\widetilde{B}}\Omega_{\widetilde{B}})=0,
\]
which implies that $d^{\widetilde{B}}\Omega_{\widetilde{B}}=0$.
\end{proof}

Note that if $G$ is the trivial group, the condition $ii)$ of
Theorem \ref{t2.4} is satisf\/ied and if only if
$\Gamma(\ker\Omega_B)$ is an ideal of $\Gamma(B).$

Finally, we deduce the following corollary.

\begin{corollary}\label{c2.5}
Let $\Psi:G\times B\to B$ be a presymplectic action of the Lie
group $G$ on the presymplectic Lie algebroid $(B,\Omega_B)$ and
suppose that:
\begin{enumerate}\itemsep=0pt
\item[i)] $\Psi_g:B\to B$ is a Lie algebroid isomorphism, for all
$g\in G.$ \item[ii)] If $X\in \Gamma(B)^p_{\widetilde{\pi}_B}$ and
$Y\in \Gamma(\ker\Omega_B)$, we have that
\[
\lcf X,Y\rcf_B\in \Gamma(\ker\Omega_B).
\]
\end{enumerate}
Then, the reduced vector bundle
$\tau_{\widetilde{B}}:\widetilde{B}\to \widetilde{N}$ admits a
unique Lie algebroid structure such that $\widetilde{\pi}_B:B\to
\widetilde{B}$ is a Lie algebroid epimorphism. Moreover, the
nondegenerate section $\Omega_{\widetilde{B}}$ induces a
symplectic structure on the Lie algebroid
$\tau_{\widetilde{B}}:\widetilde{B}\to \widetilde{N}$.
\end{corollary}
\begin{proof} If $X,Y\in \Gamma(B)^p_{\widetilde{\pi}_B}$ then
\[
\Psi_g\circ X=(X+Z_g)\circ \psi_g,\qquad \Psi_g\circ Y=(Y+W_g)\circ
\psi_g, \qquad \mbox{for all} \ \ g\in G,
\]
where $Z_g,W_g\in \Gamma(\ker\Omega_B)$ and $\psi_g:N\to N$ is the
dif\/feomorphism associated with $\Psi_g:B\to B$.  Thus, using the
fact that $\Psi_g$ is a Lie algebroid isomorphism, it follows that
\[
\Psi_g\circ \lcf X,Y\rcf_B\circ \psi_{g^{-1}}=\lcf
X+Z_g,Y+W_g\rcf_B,\qquad \mbox{for all} \ \ g\in G.
\]
Therefore, from condition $ii)$ in the corollary, we deduce that
\[
\Psi_g\circ \lcf X,Y\rcf_B\circ \psi_{g^{-1}}-\lcf X,Y\rcf_B\in
\Gamma(\ker\Omega_B),\qquad \mbox{for all} \ \ g\in G,
\]
which implies that $\lcf X,Y\rcf_B\in
\Gamma(B)^p_{\widetilde{\pi}_B}.$ This proves that
$\Gamma(B)^p_{\widetilde{\pi}_B}$ is a Lie subalgebra of the Lie
algebra $(\Gamma(B),\lcf\cdot,\cdot\rcf_B)$.
\end{proof}

\subsection{Reduction of the Hamiltonian dynamics}\label{section2.2}
Let $A$ be a Lie algebroid over a manifold $M$ and $\Omega_A$ be a
symplectic section of $A$. In addition, suppose that $B$ is a Lie
subalgebroid of~$A$ over the submanifold $N$ of $M$ and that $G$
is a Lie group such that one may construct the symplectic
reduction $\tau_{\widetilde{B}}:\widetilde{B}=(B/\ker
\Omega_B)/G\to \widetilde{N}=N/G$ of~$A$ by $B$ and $G$ as in
Section \ref{section2.1} (see Theorem \ref{t2.4}).

We will also assume that $B$ and $N$ are subsets of $A$ and $M$,
respectively (that is, the corresponding immersions $i_B:B\to A$
and $i_N:N\to M$ are the canonical inclusions), and that $N$ is a
closed submanifold.

\begin{theorem}[{\bf The reduction of the Hamiltonian dynamics}]\label{t2.6}
Let $H_M:M\to \R$ be a~Hamitonian function for the symplectic Lie
algebroid $A$ such that:
\begin{enumerate}\itemsep=0pt
\item[i)] The restriction $H_N$ of $H_M$ to $N$ is $G$-invariant
and

\item[ii)] If ${\mathcal H}_{H_M}^{\Omega_A}$ is the Hamiltonian
section of $H_M$ with respect to the symplectic section
$\Omega_A$, we have that ${\mathcal
H}_{H_M}^{\Omega_A}(N)\subseteq B$.
\end{enumerate}
Then:
\begin{enumerate}\itemsep=0pt
\item[a)] $H_N$ induces a real function
$H_{\widetilde{N}}:\widetilde{N}\to \R$ such that
$H_{\widetilde{N}}\circ \pi_N=H_N.$ \item[b)] The restriction of
${\mathcal H}_{H_M}^{\Omega_A}$ to $N$ is
$\widetilde{\pi}_B$-projectable over the Hamiltonian section of
the function $H_{\widetilde{N}}$ with respect to the reduced
symplectic structure $\Omega_{\widetilde{B}}$ and \item[c)] If
$\gamma:I\to M$ is a solution of the Hamilton equations for $H_M$
in the symplectic Lie algebroid $(A,\Omega_A)$ such that
$\gamma(t_0)\in N$, for some $t_0\in I$, then $\gamma(I)\subseteq
N$ and $\pi_N\circ \gamma:I\to \widetilde{N}$ is a~solution of the
Hamilton equations for $H_{\widetilde{N}}$ in the symplectic Lie
algebroid $(\widetilde{B},\Omega_{\widetilde{B}})$.
\end{enumerate}
\end{theorem}
\begin{proof} $a)$ Using that the function $H_N$ is $G$-invariant, we
deduce $a)$.

$b)$ If $x\in N$ and $a_x\in B_x$ then, since
$\widetilde{\pi}_B^*\Omega_{\widetilde{B}}=\Omega_B$ (see Theorem
\ref{t2.4}), ${\mathcal H}_{H_M}^{\Omega_A}(N)\subseteq B$ and
$\tau_B:B\to N$ is a Lie subalgebroid of $A$, we have that
\[
(i(\widetilde{\pi}_B({\mathcal
H}_{H_M}^{\Omega_A}(x)))\Omega_{\widetilde{B}}(\pi_N(x))
(\widetilde{\pi}_B(a_x)))=(d^BH_N)(x)(a_x).
\]
Thus, using that $\widetilde{\pi}_B$ is a Lie algebroid
epimorphism and the fact that $H_{\widetilde{N}}\circ \pi_N=H_N,$
it follows that
\begin{gather*}
(i(\widetilde{\pi}_B({\mathcal
H}_{H_M}^{\Omega_A}(x)))\Omega_{\widetilde{B}}(\pi_N(x)))
(\widetilde{\pi}_B(a_x)))=(d^{\widetilde{B}}H_{\widetilde{N}})
(\pi _N(x))(\widetilde{\pi}_B(a_x))\\
\phantom{(i(\widetilde{\pi}_B({\mathcal
H}_{H_M}^{\Omega_A}(x)))\Omega_{\widetilde{B}}(\pi_N(x)))
(\widetilde{\pi}_B(a_x)))}{}=(i({\mathcal
H}_{H_{\widetilde{N}}}^{\Omega_{\widetilde{B}}}(\pi_N(x)))
\Omega_{\widetilde{B}}(\pi_N(x)))(\widetilde{\pi}_B(a_x))).
\end{gather*}
This implies that
\[
\widetilde{\pi}_B({\mathcal H}_{H_M}^{\Omega_A}(x))={\mathcal
H}_{H_{\widetilde{N}}}^{\Omega_{\widetilde{B}}}(\pi_N(x)).
\]

$c)$ Using $b)$, we deduce that the vector f\/ield
$\rho_B(({\mathcal H}_{H_M}^{\Omega_A})_{|N})$ is
$\pi_N$-projectable on the vector f\/ield
$\rho_{\widetilde{B}}({\mathcal
H}_{H_{\widetilde{N}}}^{\Omega_{\widetilde{B}}})$ (it is a
consequence of the equality $\rho_{\widetilde{B}}\circ
\widetilde{\pi}_B=T\pi_N\circ \rho_B$). This proves $c).$ Note
that $N$ is closed and that the integral curves of
$\rho_B(({\mathcal H}_{H_M}^{\Omega_A})_{|N})$ (respectively,
$\rho_{\widetilde{B}}({\mathcal
H}_{H_{\widetilde{N}}}^{\Omega_{\widetilde{B}}})$) are the
solutions of the Hamilton equations for $H_M$ in $A$
(respectively, for $H_{\widetilde{N}}$ in $\widetilde{B}$) with
initial condition in $N$ (respectively, in $\widetilde{N}$).
\end{proof}

\section{Examples and applications}\label{section3}

\subsection{Cartan  symplectic reduction in the presence of a symmetry Lie group}\label{section3.1}

Let $M$ be a symplectic manifold with symplectic $2$-form $\Omega_{TM}$. Suppose that $N$ is a submanifold of $M$, that $G$ is a
Lie group and that $N$ is the total space of a principal $G$-bundle over $\widetilde{N}=N/G.$ We will denote $\psi:G\times N\to N$
the free action of $G$ on $N$, by $\pi_N:N\to \widetilde{N}=N/G$ the principal bundle projection and by $\Omega_{TN}$ the $2$-form
on $N$ given by
\[
\Omega_{TN}=i_N^*\Omega_{TM}.
\]
Here, $i_N:N\to M$ is the canonical inclusion.

\begin{proposition}\label{p3.1}
If the vertical bundle to $\pi_N$ is the kernel of the $2$-form
$\Omega_{TN}$, that is,
\begin{equation}
\label{3.1}
V\pi_N=\ker\Omega_{TN},
\end{equation}
then there exists a unique symplectic $2$-form $\Omega_{T(N/G)}$ on
$\widetilde{N}=N/G$ such that $\pi_N^*\Omega_{T(N/G)}=\Omega_{TN}.$
\end{proposition}

\begin{proof} We will prove the result using Theorem~\ref{t2.3} and
Corollary~\ref{c2.5}. In fact, we apply Theorem~\ref{t2.3} to the
followings elements:
\begin{itemize}\itemsep=0pt
\item The standard symplectic Lie algebroid $\tau_{TM}:TM\to M,$
\item The Lie subalgebroid of $TM$, $\tau_{TN}:TN\to N$ and \item
The presymplectic action $T\psi:G\times TN\to TN$, i.e., the
tangent lift of the action $\psi$ of $G$ on $N$.
\end{itemize}

The corresponding reduced vector bundle
$(\widetilde{TN}=(TN/\ker\Omega_{TN})/G, \Omega_{\widetilde{TN}})$
is isomorphic to the tangent bundle $T(N/G)$ in a natural way.

On the other hand, it is well-known that the Lie bracket of a
$\pi_N$-projectable vector f\/ield and a $\pi_N$-vertical vector
f\/ield is a $\pi_N$-vertical vector f\/ield. Therefore, we may use
Corollary \ref{c2.5} and we deduce that the symplectic vector
bundle $\tau_{\widetilde{TN}}: \widetilde{TN}\to\widetilde{N}=N/G$
admits a unique Lie algebroid structure such that the canonical
projection  $\widetilde{\pi}_{TN}:TN\to \widetilde{TN}$ is a Lie
algebroid epimorphism. In addition, $\Omega_{\widetilde{TN}}$ is a
symplectic section of the Lie algebroid $\tau_{\widetilde{TN}}:
\widetilde{TN}\to\widetilde{N}=N/G$.

Finally, under the identif\/ication between the vector bundles
$\widetilde{TN}$ and $T(N/G)$, one may see that the resultant Lie
algebroid structure on the vector bundle $\tau_{T(N/G)}: T(N/G)\to
N/G$ is the  just standard Lie algebroid structure.
\end{proof}

\begin{remark}\label{r3.2}
 The reduction process described in Proposition \ref{p3.1} is
just the classical Cartan symplectic reduction process (see, for
example, \cite{OR}) for the particular case when the kernel of the
presymplectic $2$-form on the submanifold $N$ of the original
symplectic manifold is just the union of the tangent spaces of the
$G$-orbits associated with a principal $G$-bundle with total space
$N.$
\end{remark}

Now, assume that the submanifold $N$ is closed. Then, using Theorem
\ref{t2.6}, we may prove the following result:

\begin{corollary}
Let $H_M:M\to \R$ be a Hamiltonian function on the symplectic
manifold $M$ such that:
\begin{enumerate}
\itemsep=0pt
\item The restriction $H_N$ of $H_M$ to $N$ is $G$-invariant and
\item The restriction to $N$  of the Hamiltonian vector field
${\mathcal H}_{H_M}^{\Omega_{TM}}$ of $H_M$ with respect to
$\Omega_{TM}$ is tangent to $N.$
\end{enumerate}
Then:
\begin{enumerate}\itemsep=0pt
\item[a)] $H_N$ induces a real function $H_{\widetilde{N}}:\widetilde{N}=N/G\to \R$ such that $H_{\widetilde{N}}
\circ \pi_N=H_N;$
\item[b)] The vector field $({\mathcal H}_{H_M}^{\Omega_{TM}})_{|N}$ on $N$ is $\pi_N$-projectable
on the Hamiltonian vector field of $H_{\widetilde{N}}$ with respect
to the reduced symplectic $2$-form $\Omega_{T(N/G)}$ on $N/G$ and
\item[c)]
If $\gamma:I\to M$ is a solution of the Hamiltonian equations for
$H_M$ in the symplectic manifold $(M,\Omega_{TM})$ such that
$\gamma(t_0)\in N$, for some $t_0\in I$, then $\gamma(I)\subseteq N$
and $\pi_N\circ \gamma:I\to \widetilde{N}$ is a~solution of the
Hamiltonian equations for $H_{\widetilde{N}}$ in the symplectic
manifold $(\widetilde{N},\Omega_{\widetilde{TN}}).$
  \end{enumerate}

\end{corollary}

\subsection{Symplectic reduction of symplectic Lie algebroids by Lie subalgebroids}\label{section3.2}

Let $A$ be a symplectic Lie algebroid over the manifold $M$ with
symplectic section $\Omega_A$ and $B$ be a Lie subalgebroid of $A$
over the submanifold $N$ of $M.$

Denote by $(\lcf\cdot,\cdot\rcf_B,\rho_B)$ the Lie algebroid
structure of $B$ and by $\Omega_B$ the presymplectic section on
$B$ given by{\samepage
\[
\Omega_B=i_B^*\Omega_A,
\]
where $i_B:B\to A$ is the canonical inclusion.}

As in Section $2$, we will assume that the dimension of the
subspace $\ker\Omega_B(x)$ is constant, for all $x\in N.$ Then,
one may consider the vector bundle
$\tau_{\ker\Omega_B}:\ker\Omega_B\to N$ which is a Lie
subalgebroid of the Lie algebroid $\tau_B:B\to N.$

\begin{corollary}\label{c3.4}
If the space $\Gamma(\ker\Omega_B)$ is an ideal of the Lie algebra
$(\Gamma(B),\lcf\cdot,\cdot\rcf_B)$ then the quotient vector
bundle $\tau_{B/\ker\Omega_B}:\widetilde{B}=B/\ker\Omega_B\to N$
admits a unique Lie algebroid structure such that the canonical
projection $\pi_B:B\to \widetilde{B}=B/\ker\Omega_B$ is a Lie
algebroid epimorphism over the identity $Id:N\to N$ of $N$.
Moreover, there exists a unique symplectic section
$\Omega_{\widetilde{B}}$ on the Lie algebroid
$\tau_{\widetilde{B}}:\widetilde{B}\to N$ which satisfies the
condition
\[
\pi_B^*\Omega_{\widetilde{B}}=\Omega_B.
\]
\end{corollary}
\begin{proof}
 It follows using Theorem \ref{t2.4} (in this case, $G$ is the trivial Lie group $G=\{e\}$).
\end{proof}

Next, we will consider the particular case when the manifold $M$ is a single point.

\begin{corollary}\label{c3.5}
Let $({\mathfrak g},[\cdot,\cdot]_{\mathfrak g})$ be a symplectic
Lie algebra of finite dimension with symplectic $2$-form
$\Omega_{\mathfrak g}:{\mathfrak g}\times {\mathfrak g}\to \R$. If
${\mathfrak h}$ is a Lie subalgebra of ${\mathfrak g}$,
$\Omega_{\mathfrak h}:{\mathfrak h}\times {\mathfrak h}\to \R$ is
the restriction of $\Omega_{\mathfrak g}$ to ${\mathfrak h}\times
{\mathfrak h}$ and $\ker\Omega_{\mathfrak h}$ is an ideal of the Lie
algebra $({\mathfrak h},[\cdot,\cdot]_{\mathfrak h})$ then:
\begin{enumerate}\itemsep=0pt
\item The quotient vector space $\widetilde{\mathfrak h}={\mathfrak h}/\ker\Omega_{\mathfrak h}$
admits a unique Lie algebra structure such that the canonical
projection $\pi_{\mathfrak h}:{\mathfrak h}\to \widetilde{\mathfrak
h}= {\mathfrak h}/\ker\Omega_{\mathfrak h}$ is a Lie algebra
epimorphism and
\item There exists a unique symplectic $2$-form $\Omega_{\widetilde{\mathfrak
h}}:\widetilde{\mathfrak h}\times \widetilde{\mathfrak h}\to \R$ on
$\widetilde{\mathfrak h}$ which satisfies the condition
\[
\pi_{\mathfrak h}^*\Omega_{\widetilde{\mathfrak
h}}=\Omega_{\mathfrak h}.
\]
\end{enumerate}
\end{corollary}
\begin{proof}
From Corollary \ref{c3.4} we deduce the result.
\end{proof}

The symplectic Lie algebra $(\widetilde{\mathfrak
h},\Omega_{\widetilde{\mathfrak h}})$ is called the \emph{symplectic
reduction of the Lie algebra $({\mathfrak g},\Omega_{\mathfrak g})$
by the Lie subalgebra ${\mathfrak h}$}.

This reduction process of symplectic Lie algebras plays an important
role in the description of symplectic Lie groups (see, for instance,
\cite{Medina}).

\subsection{Reduction of a symplectic Lie algebroid by a symplectic Lie subalgebroid\\ and a symmetry Lie group}\label{section3.3}

Let $A$ be a symplectic Lie algebroid over the manifold $M$ with
symplectic section $\Omega_A$ and $B$ be a Lie subalgebroid of $A$
over the submanifold $N$ of $M.$

Denote by $(\lcf\cdot,\cdot\rcf_B,\rho_B)$ the Lie algebroid
structure on $B$ and by $\Omega_B$ the presymplectic section on
$B$ given by
\[
\Omega_B=i_B^*\Omega_A,
\]
where $i_B:B\to A$ is the canonical inclusion.

We will assume that $\Omega_B$ is nondegenerate. Thus, $\Omega_B$
is a symplectic section on the Lie subalgebroid $\tau_B:B\to N$
and $\ker\Omega_B(x)=\{0\}$, for all $x\in N.$

\begin{corollary}\label{c3.6}
Let $\Psi:G\times B\to B$ be a symplectic action of the Lie group
$G$ on the symplectic Lie algebroid $(B,\Omega_B)$ such that the
space $\Gamma(B)^G$ of $G$-invariant sections of $B$,
\[
\Gamma(B)^G=\{X\in \Gamma(B)/\Psi_g\circ X=X\circ \psi_g,\mbox{
for all }g\in G\}
\]
is a Lie subalgebra of $(\Gamma(B),\lcf\cdot,\cdot\rcf_B)$, where
$\psi:G\times N\to N$ the corresponding action on $N.$ Then:
\begin{enumerate}
\item The quotient vector bundle
$\tau_{\widetilde{B}}:\widetilde{B}=B/G\to \widetilde{N}=N/G$
admits a unique Lie algebroid structure such that the canonical
projection $\widetilde{\pi}_B:B\to \widetilde{B}$ is a Lie
algebroid epimorphism and \item There exists a unique symplectic
structure $\Omega_{\widetilde{B}}$ on the Lie algebroid
$\tau_{\widetilde{B}}:\widetilde{B}\to \widetilde{N}$ such that
\[
\widetilde{\pi}_B^*\Omega_{\widetilde{B}}=\Omega_B.
\]
\end{enumerate}
\end{corollary}

\begin{proof}It follows using Theorem \ref{t2.4}.
\end{proof}

\begin{remark}\label{r3.6'}{\rm Let $A$ be a Lie algebroid over a manifold $M$ with Lie algebroid structure $(\lcf\cdot,\cdot\rcf_A,
\rho_A)$ and vector bundle projection $\tau_A:A\to M.$

Suppose that we have a principal action $\psi:G\times M\to M$ of
the Lie group $G$ on $M$ and an action $\Psi:G\times A\to A$ of
$G$ on $A$ such that $\tau_A:A\to M$ is $G$-equivariant and $\Psi$
is a Lie algebroid action.

Then, we may consider the quotient vector bundle $A/G$ over $M/G$.
In fact, the canonical projection $\pi_A:A\to A/G$ is a vector
bundle morphism and $(\pi_A)_{|A_x}:A_x\to (A/G)_{\pi_M(x)}$ is a~linear isomorphism, for all $x\in M.$ As we know, the space
$\Gamma(A/G)$ of sections of the vector bundle $\tau_{A/G}:A/G\to
M/G$  may be identif\/ied with the set $\Gamma(A)^G$ of
$G$-invariant sections of $A$, that is,
\[
\Gamma(A/G)\cong \Gamma(A)^G=\{X\in \Gamma(A)/\Psi_g\circ X=X\circ \psi_g,\ \forall\, g\in G\}.
\]
Now, if $g\in G$ and $X,Y\in \Gamma(A)^G$ then, since $\Psi_g$ is
a Lie algebroid morphism, we deduce that
\[
\Psi_g\circ \lcf X,Y\rcf_A\circ \psi_{g^{-1}}=\lcf\Psi_g\circ
X\circ \psi_{g^{-1}},\Psi_g\circ Y\circ \psi_{g^{-1}}\rcf_A=\lcf
X,Y\rcf_A,\]
 i.e., $\lcf X,Y\rcf_A\in \Gamma(A)^G.$ Thus,
$\Gamma(A)^G$ is a Lie subalgebra of
$(\Gamma(A),\lcf\cdot,\cdot\rcf_A)$ and therefore, using Theorem
\ref{tA1} (see Appendix), it follows that the quotient vector
bundle $A/G$ admits a unique Lie algebroid structure such that
$\pi_A$ is a Lie algebroid morphism.}
\end{remark}

Next, we will apply Corollary \ref{c3.6} to the following example.

\begin{example}\label{e3.6''}
{\rm Let $A$ be a Lie algebroid over a manifold $M$ and $G$ be a
Lie group as in Remark~\ref{r3.6'}. Note that the action $\Psi$ of
$G$ on $A$ is principal and it induces an action $\Psi^*$ of $G$
on $A^*$ which is also principal. In fact,
$\Psi_g^*=\Psi_{g^{-1}}^t,$ where $\Psi_{g^{-1}}^t$ is the dual
map of $\Psi_{g^{-1}}:A\to A$. The space $A^*/G$ of orbits of this
action is a quotient vector bundle over $M/G$ which is isomorphic
to the dual vector bundle to $\tau_{A/G}:A/G\to M/G.$ Moreover,
the canonical projection $\pi_{A^*}:A^*\to A^*/G$ is a vector
bundle morphism and $(\pi_{A^*})_{|A^*_x}:A_x^*\to
(A^*/G)_{\pi_M(x)}$ is a linear isomorphism, for all $x\in M.$

Now, we consider the action ${\mathcal T}^*\Psi$ of $G$ on the
$A$-tangent bundle to $A^*$, ${\mathcal T}^AA^*$, given by
\[
{\mathcal T}^*\Psi(g,(a_x,X_{\alpha_x}))=(\Psi_g(a_x),
(T_{\alpha_x}\Psi_{g^{-1}}^t)(X_{\alpha_x})),
\]
for $g\in G$ and $(a_x,X_{\alpha_x})\in {\mathcal
T}^A_{\alpha_x}A^*.$ Since $\Psi_g$ and $T\Psi_{g^{-1}}^t$ are Lie
algebroid isomorphisms, we have that ${\mathcal T}_g^*\Psi$ is
also a Lie algebroid isomorphism. Therefore, from Remark
\ref{r3.6'}, we have that the space
\[
\Gamma({\mathcal T}^AA^*)^G=\{\widetilde{X}\in \Gamma({\mathcal
T}^AA^*)/{\mathcal T}^*_g\Psi\circ \widetilde{X}=\widetilde{X}
\circ \Psi_{g^{-1}}^t,\ \forall\, g\in G\}
\]
of $G$-invariant sections of $\tau_{{\mathcal T}^AA^*}:{\mathcal
T}^AA^*\to A^*$ is a Lie subalgebra of the Lie algebra
$(\Gamma({\mathcal T}^AA^* ),$ $\lcf\cdot,\cdot\rcf_{{\mathcal
T}^AA^*}).$

In addition, if $\Theta_{{\mathcal T}^AA^*}\in \Gamma(({\mathcal
T}^AA^*)^*)$ is the Liouville section associated with $A$, it is
easy to prove that
\[
({\mathcal T}_g^*\Psi)^*(\Theta_{{\mathcal
T}^AA^*})=\Theta_{{\mathcal T}^AA^*},\qquad \mbox{for all}\ \ g\in G.
\]
Thus, if $\Omega_{{\mathcal T}^AA^*}\in \Gamma(\wedge^2({\mathcal
T}^AA^*)^*)$ is the canonical symplectic section associated with
$A$, it follows that
\[
({\mathcal T}^*_g\Psi)^*(\Omega_{{\mathcal
T}^AA^*})=\Omega_{{\mathcal T}^AA^*},\qquad \mbox{for all} \ \ g\in
G.
\]
This implies that ${\mathcal T}^*\Psi$ is a symplectic action of
$G$ on ${\mathcal T}^AA^*.$ Consequently, using Corollary~\ref{c3.6}, we deduce that the quotient vector bundle
$\tau_{{\mathcal T}^AA^*/G}:{\mathcal T}^AA^*/G\to A^*/G$ admits a
unique Lie algebroid structure such that the canonical projection
$\pi_{{\mathcal T}^AA^*}:{\mathcal T}^AA^*\to {\mathcal T}^AA^*/G$
is a Lie algebroid epimorphism. Moreover, there exists a unique
symplectic structure $\Omega_{({\mathcal T}^AA^*)/G}$ on the Lie
algebroid $\tau_{({\mathcal T}^AA^*)/G}:({\mathcal T}^AA^*)/G\to
A^*/G$ such that
\[
\pi_{{\mathcal T}^AA^*}^*(\Omega_{({\mathcal
T}^AA^*)/G})=\Omega_{{\mathcal T}^AA^*}.
\]
Next, we will prove that the symplectic Lie algebroid $({\mathcal
T}^AA^*/G,\Omega_{({\mathcal T}^AA^*)/G})$ is isomorphic to the
$A/G$-tangent bundle to $A^*/G$. For this purpose, we will
introduce the map $(\pi_A,T\pi_{A^*}):{\mathcal T}^AA^*\to
{\mathcal T}^{A/G} (A^*/G)$ def\/ined by
\[
(\pi_A, T\pi_{A^*})(a_x,X_{\alpha_x})=(\pi_A(a_x),
(T_{\alpha_x}\pi_{A^*})(X_{\alpha_x})), \qquad \mbox{for} \ \
(a_x,X_{\alpha_x})\in {\mathcal T}_{\alpha_x}^AA^*.
\]
We have that $(\pi_A,T\pi_{A^*})$ is a vector bundle morphism.
Furthermore, if $(a_x,X_{\alpha_x})\kern-2pt\in\kern-1pt {\mathcal
T}_{\alpha_x}^AA^*$ and
\[
0=(\pi_A,T\pi_{A^*})(a_x,X_{\alpha_x})
\]
then $a_x=0$ and, therefore, $X_{\alpha_x}\in T_{\alpha_x}A^*_x$
and $(T_{\alpha_x}\pi_{A^*})(X_{\alpha_x})=0$. But, as
$(\pi_{A^*})_{|A_x^*}:A_x^*\to (A^*/G)_{\pi_M(x)}$ is a linear
isomorphism, we conclude that $X_{\alpha_x}=0.$ Thus,
$(\pi_A,T\pi_{A^*})_{|{\mathcal T}_{\alpha_x}^AA^*}:{\mathcal
T}_{\alpha_x}^AA^*\to {\mathcal
T}_{\pi_{A^*}(\alpha_x)}^{A/G}(A^*/G)$ is a linear isomorphism
(note that $\dim {\mathcal T}_{\alpha_x}^AA^*=\dim {\mathcal
T}_{\pi_{A^*}(\alpha)}^{A/G}(A^*/G))$.

On the other hand, using that $\pi_A$ (respectively, $T\pi_{A^*}$)
is a Lie algebroid morphism between the Lie algebroids
$\tau_A:A\to M$ and $\tau_{A/G}:A/G\to M/G$ (respectively,
$\tau_{TA^*}:TA^*\to A^*$ and $\tau_{T(A^*/G)}:T(A^*/G)\to A^*/G$)
we deduce that $(\pi_A,T\pi_{A^*})$ is also a Lie algebroid
morphism $M.$

Moreover, if $\Theta_{{\mathcal T}^{A/G}(A^*/G)}\in
\Gamma({\mathcal T}^{A/G}(A^*/G))^*$ is the Liouville section
associated with the Lie algebroid $\tau_{A/G}:A/G\to M/G$, we have
that
\[
(\pi_A,T\pi_{A^*})^*(\Theta_{{\mathcal
T}^{A/G}(A^*/G)})=\Theta_{{\mathcal T}^AA^*}
\]
which implies that
\[
(\pi_A,T\pi_{A^*})^*(\Omega_{{\mathcal
T}^{A/G}(A^*/G)})=\Omega_{{\mathcal T}^AA^*},
\]
where $\Omega_{{\mathcal T}^{A/G}(A^*/G)}$ is the canonical
symplectic section associated with the Lie algebroid
$\tau_{A/G}:A/G\to M/G$.

Now, since $\pi_A\circ \Psi_g=\pi_A$ and $\pi_{A^*}\circ
\Psi_g^*=\pi_{A^*}$ for all $g\in G,$ we deduce that the map
$(\pi_A,T\pi_{A^*}): {\mathcal T}^AA^*\to {\mathcal
T}^{A/G}{(A^*/G)}$ induces a map
$(\widetilde{\pi_A,T\pi_{A^*}}):({\mathcal T}^AA^*)/G\to {\mathcal
T}^{A/G}(A^*/G)$ such that
\[
\widetilde{(\pi_A,T\pi_{A^*})}\circ \pi_{{\mathcal
T}^AA^*}=(\pi_A,T\pi_{A^*}).
\]
Finally, using the above results, we obtain
 that $\widetilde{(\pi_A,T\pi_{A^*})}$ is a Lie algebroid isomorphism and, in addition,
\[
\widetilde{(\pi_A,T\pi_{A^*})}^*\Omega_{{\mathcal
T}^{A/G}(A^*/G)}=\Omega_{({\mathcal T}^AA^*)/G}.
\]
In conclusion, we have proved that the symplectic Lie algebroids
$(({\mathcal T}^AA^*)/G,\Omega_{({\mathcal T}^AA^*)/G})$ and
$({\mathcal T}^{A/G}(A^*/G),\Omega_{{\mathcal T}^{A/G}(A^*/G)})$
are isomorphic. }\end{example} Next, we will prove the following
result.

\begin{corollary}\label{c3.7}
Under the same hypotheses as in Corollary {\rm \ref{c3.4}}, if the
submanifold $N$ is closed and $H_M:M\to \R$ is a Hamiltonian
function for the symplectic Lie algebroid $(A,\Omega_A)$ such that
the restriction $H_N$ of $H_M$ is $G$-invariant then:

\begin{enumerate}\itemsep=0pt
\item $H_N$ induces a real function
$H_{\widetilde{N}}:\widetilde{N}=N/G\to \R$ such that
$H_{\widetilde{N}}\circ \pi_N=H_N;$ \item The Hamiltonian section
${\mathcal H}_{H_M}^{\Omega_A}$ of $H_M$ in the symplectic Lie
algebroid $(A,\Omega_A)$ satisfies the condition
\[
{\mathcal H}_{H_M}^{\Omega_A}(N)\subseteq B.
\]
Moreover, $({\mathcal H}_{H_M}^{\Omega_A})_{|N}$ is
$\widetilde{\pi}_B$-projectable on the Hamiltonian section of the
function $H_{\widetilde{N}}$ with respect to the reduced
symplectic section $\Omega_{\widetilde{B}}$ and \item If
$\gamma:I\to M$ is a solution of the  Hamilton equations for $H_M$
in the symplectic Lie algebroid $(A,\Omega_A)$ such that
$\gamma(t_0)\in N$  for some $t_0\in I$, then $\gamma(I)\subseteq
N$ and $\pi_N\circ \gamma:I\to \widetilde{N}$ is a~solution of the
Hamilton equations for $H_{\widetilde{N}}$ in the symplectic Lie
algebroid $(\widetilde{B},\Omega_{\widetilde{B}})$.
\end{enumerate}
\end{corollary}

\begin{proof} Let ${\mathcal H}_{H_N}^{\Omega_B}$ be the Hamiltonian
section of the Hamiltonian function $H_N:N\to \R$  in the
symplectic Lie algebroid $(B,\Omega_B)$. Then, we have that
$({\mathcal H}_{H_M}^{\Omega_A})_{|N}{=}{\mathcal H}_N^{\Omega_B}$
and, thus, ${\mathcal H}_{H_M}^{\Omega_A}(N){\subseteq} B.$
Therefore, using Theorem \ref{t2.6}, we deduce the result.
\end{proof}

Next, we will apply the above results in order to give a ``symplectic
description'' of the Hamiltonian reduction process by stages in the
Poisson setting.

This reduction process may be described as follows (see, for
instance, \cite{CMR,OR}).

Let $Q$ be the total space of a principal $G_1$-bundle and the
conf\/iguration space of an standard Hamiltonian system with
Hamiltonian function $H:T^*Q\to \R$. We will assume that $H$ is
$G_1$-invariant. Then, the quotient manifold $T^*Q/G_1$ is a Poisson
manifold (see Section \ref{section1.2}) and if $\pi_{T^*Q}^1:T^*Q\to
T^*Q/G_1$ is the canonical projection, we have that the Hamiltonian
system $(T^*Q,H)$ is $\pi_{T^*Q}^1$-projectable over a Hamiltonian
system on $T^*Q/G_1$. This means that $\pi_{T^*Q}^1:T^*Q\to
T^*Q/G_1$ is a Poisson morphism (see Section \ref{section1.2}) and
that there exists a real function $\widetilde{H}_1:T^*Q/G_1\to \R$
such that $\widetilde{H}_1\circ \pi_{T^*Q}^1=H.$ Thus, if
$\gamma:I\to T^*Q$ is a solution of the Hamilton equations for $H$
in $T^*Q$ then $\pi_{T^*Q}^1\circ \gamma:I\to T^*Q/G_1$ is a
solution of the Hamilton equations for $\widetilde{H}_1$ in
$T^*Q/G_1.$

Now, suppose that $G_2$ is a closed normal subgroup of $G_1$ such
that the action of $G_2$ on $Q$ induces a principal $G_2$-bundle.
Then, it is clear that the construction of the reduced Hamiltonian
system $(T^*Q/G_1,\widetilde{H}_1)$ can be carried out in two steps.

{\bfseries \itshape  First step:} The Hamiltonian function $H$ is
$G_2$-invariant and, therefore, the Hamiltonian system $(T^*Q,H)$
may be reduced to a Hamiltonian system in the Poisson manifold
$T^*Q/G_2$ with Hamiltonian function $\widetilde{H}_2:T^*Q/G_2\to
\R$.

{\bfseries \itshape Second step:} The cotangent lift of the action of $G_1$
on $Q$ induces a Poisson action of the quotient Lie group
$G_1/G_2$ on the Poisson manifold $T^*Q/G_2$ and the space
$(T^*Q/G_2)/(G_1/G_2)$ of orbits of this action is isomorphic to
the full reduced space $T^*Q/G_1$. Moreover, the function~$\widetilde{H}_2$ is $(G_1/G_2)$-invariant and, consequently, it
induces a Hamiltonian function $\widetilde{\widetilde{H}}_2$ on
the space $(T^*Q/G_2)/(G_1/G_2).$ Under the identif\/ication between
$(T^*Q/G_2)/(G_1/G_2)$ and $T^*Q/G_1$,
$\widetilde{\widetilde{H}}_2$ is just the Hamiltonian function
$\widetilde{H}_1$. In conclusion, in this second step, we see that
the Hamiltonian system $(T^*Q/G_2,\widetilde{H}_2)$ may be reduced
to the Hamiltonian system $(T^*Q/G_1,\widetilde{H}_1).$

Next, we will give a ``symplectic description" of the above
reduction process.

For this purpose, we will consider the Atiyah algebroids
$\tau_{TQ/G_2}:TQ/G_2\to \widetilde{Q}_2=Q/G_2$ and
$\tau_{TQ/G_1}:TQ/G_1\to \widetilde{Q}_1=Q/G_1$ associated with the
principal bundles $\pi^2_Q:Q\to \widetilde{Q}_2=Q/G_2$ and
$\pi^1_Q:Q\to \widetilde{Q}_1=Q/G_1$. We also consider the
corresponding dual vector bundles $\tau_{T^*Q/G_2}:T^*Q/G_2\to
\widetilde{Q}_2=Q/G_2$ and $\tau_{T^*Q/G_1}:T^*Q/G_1\to
\widetilde{Q}_1=Q/G_1$.

Then, the original Hamiltonian system $(T^*Q,H)$ may be considered
as a Hamiltonian system, with Hamiltonian function $H$, in the
standard symplectic Lie algebroid
\[
\tau_{T(T^*Q)}:T(T^*Q)\to T^*Q.
\]
Note that this symplectic Lie algebroid is isomorphic to the
$TQ$-tangent bundle to $T^*Q.$

On the other hand, the reduced Poisson Hamiltonian system
$(T^*Q/G_1,\widetilde{H}_1)$ may be considered as a Hamiltonian
system, with Hamiltonian function $\widetilde{H}_1:T^*Q/G_1\to \R$
in the symplectic Lie algebroid
\[
\tau_{{\mathcal T}^{TQ/G_1}(T^*Q/G_1)}:{\mathcal T}^{TQ/G_1}(T^*Q/G_1)\to T^*Q/G_1.
\]
Now, in order to give the symplectic description of the two steps of
the reduction process, we proceed as follows:

{\bfseries \itshape  First step:} The action of $G_2$ on $Q$ induces a
symplectic action of $G_2$ on the symplectic Lie algebroid
$\tau_{{\mathcal T}^{TQ}(T^*Q)}: {\mathcal T}^{TQ}(T^*Q)\cong
T(T^*Q)\to T^*Q.$ Furthermore (see Example \ref{e3.6''}), the
symplectic reduction of this Lie algebroid by $G_2$ is just the
$TQ/G_2$-tangent bundle to $T^*Q/G_2$, that is, ${\mathcal
T}^{TQ/G_2}(T^*Q/G_2)$ (note that the standard Lie bracket of two
$G_2$-invariant vector f\/ields on~$Q$ is again a $G_2$-invariant
vector f\/ield). Therefore, the Hamiltonian system, in the symplectic
Lie algebroid $\tau_{T(T^*Q)}:T(T^*Q)\to T^*Q,$ with Hamiltonian
function $H:T^*Q\to \R$ may be reduced to a Hamiltonian system, in
the symplectic Lie algebroid $\tau_{{\mathcal
T}^{TQ/G_2}(T^*Q/G_2)}:{\mathcal T}^{TQ/G_2}(T^*Q/G_2)\to T^*Q/G_2,$
with Hamiltonian function $\widetilde{H}_2:T^*Q/G_2\to \R.$

{\bfseries \itshape  Second step:} The action of the Lie group $G_1$ on $Q$
induces a Lie algebroid action of the quotient Lie group $G_1/G_2$
on the Atiyah algebroid $\tau_{TQ/G_2}:TQ/G_2\to
\widetilde{Q}_2=Q/G_2$. Thus (see Remark \ref{r3.6'}), the
quotient vector bundle $(TQ/G_2)/(G_1/G_2)\to (Q/G_2)/(G_1/G_2)$
admits a quotient Lie algebroid structure. If follows that this
Lie algebroid is isomorphic to the Atiyah algebroid
$\tau_{TQ/G_1}:TQ/G_1\to \widetilde{Q}_1=Q/G_1.$ On the other
hand, the Lie algebroid action of $G_1/G_2$ on the Atiyah
algebroid $\tau_{TQ/G_2}:TQ/G_2\to \widetilde{Q}_2=Q/G_2$ induces
a~symplectic action of $G_1/G_2$ on the symplectic Lie algebroid
${\mathcal T}^{TQ/G_2}(T^*Q/G_2)\to T^*Q/G_2$ which is also a Lie
algebroid action (see Example \ref{e3.6''}). Moreover (see again
Example \ref{e3.6''}), the symplectic reduction of
$\tau_{{\mathcal T}^{TQ/G_2}(T^*Q/G_2)}:{\mathcal
T}^{TQ/G_2}(T^*Q/G_2)\to T^*Q/G_2$ by $G_1/G_2$ is the
$(TQ/G_2)/(G_1/G_2)$-tangent bundle to $(T^*Q/G_2)/(G_1/G_2)$,
that is, the symplectic Lie algebroid $\tau_{{\mathcal
T}^{TQ/G_1}(T^*Q/G_1)}:{\mathcal T}^{TQ/G_1}(T^*Q/G_1)\to
T^*Q/G_1$. Therefore, the Hamiltonian system in the symplectic Lie
algebroid $\tau_{{\mathcal T}^{TQ/G_2}(T^*Q/G_2)}:{\mathcal
T}^{TQ/G_2}(T^*Q/G_2)\to T^*Q/G_2$, with Hamiltonian function
$\widetilde{H}_2:T^*Q/G_2\to \R$, may be reduced to the
Hamiltonian system in the symplectic Lie algebroid
$\tau_{{\mathcal T}^{TQ/G_1}(T^*Q/G_1)}:{\mathcal
T}^{TQ/G_1}(T^*Q/G_1)\to T^*Q/G_1$, with Hamiltonian function
$\widetilde{H}_1:T^*Q/G_1\to \R.$

\subsection{ A Lagrange top}\label{section3.4}

In Sections \ref{section3.2} and \ref{section3.3}, we discussed
the reduction of a symplectic Lie algebroid $A$ by a Lie
subalgebroid $B$ and a symmetry Lie group $G$ for the particular
case when $G$ is the trivial group and for the particular case
when $B$ is symplectic, respectively.

Next, we will consider an example such that the Lie subalgebroid
$B$ is not symplectic and the symmetry Lie group $G$ is not
trivial: a Lagrange top.

\emph{A Lagrange top } is a classical mechanical system which
consists of a symmetric rigid body with a f\/ixed point moving in a
gravitational f\/ield (see, for instance, \cite{BS,Marsden}).

A Lagrange top may be described as a Hamiltonian system on the
$\widetilde{A}$-tangent bundle to $\widetilde{A}^*,$ where
$\widetilde{A}$ is an action Lie algebroid over the sphere $S^2$
in $\R^3$
\[
S^2=\{\vec{x}\in \R^3/\|\vec{x}\|=1\}.
\]
The Lie algebroid structure on $\widetilde{A}$ is def\/ined as
follows.

Let $SO(3)$ be the special orthogonal group of order $3$. As we
know, the Lie algebra ${\mathfrak{so}}(3)$ of $SO(3)$ may identif\/ied
with $\R^3$ and, under this identif\/ication, the Lie bracket
$[\cdot,\cdot]_{{\mathfrak{so}}(3)}$ is just the cross product on
$\R^3.$

We have the standard left action of $SO(3)$ on $S^2$ and the
corresponding inf\/initesimal left action $\Phi:{\mathfrak{so}}(3)
\cong \R^3\to {\mathfrak X}(S^2)$ given by
\[
\Phi(\xi)(\vec{x})=\xi\times \vec{x},\qquad \mbox{for} \ \ \xi \in
\mathfrak{so}(3)\cong \R^3 \ \ \mbox{and} \ \ \vec{x}\in S^2.
\]
$\widetilde{A}$ is the Lie algebroid over $S^2$ associated with
the inf\/initesimal left action $\Phi$. Thus, $\widetilde{A}=
S^2\times {\mathfrak{so}}(3)\cong S^2\times \R^3$ and if
$\lcf\cdot,\cdot\rcf_{\widetilde{A}}$  is the Lie bracket on the
space $\Gamma(\widetilde{A})$ we can choose a global basis
$\{e_1,e_2,e_3\}$ of $\Gamma(\widetilde{A})$ such that
\[
\lcf e_1,e_2\rcf_{\widetilde{A}}=e_3,\qquad \lcf
e_3,e_1\rcf_{\widetilde A}=e_2,\qquad \lcf
e_2,e_3\rcf_{\widetilde{A}}=e_1.
\]
The dual vector bundle $\widetilde{A}^*$ to $\widetilde{A}$ is the
trivial vector bundle $\tau_{\widetilde{A}^*}:\widetilde{A}^*
\cong S^2\times \R^3\to S^2$ and the Hamiltonian function
$\widetilde{H}:\widetilde{A}^*\cong S^2\times \R^3\to \R$ is
def\/ined by
\[
\widetilde{H}(\vec{x},\pi)=\frac{1}{2}\left(\frac{\pi_1^2}{I} + \frac{\pi_2^2}{I} + \frac{\pi_3^2}{J}\right) + (mgl)z,
\]
for $\vec{x}=(x,y,z)\in S^2$ and $\pi=(\pi_1,\pi_2,\pi_3)\in \R^3,$ where $(I,I,J)$ are the (positive) eigenvalues of the inertia
tensor, m is the mass, $g$ is the gravity and $l$ is the distance from the f\/ixed point of the body to the centre of mass.

The $\widetilde{A}$-tangent bundle to $\widetilde{A}^*,$
${\mathcal T}^{\widetilde{A}}\widetilde{A}^*,$ is isomorphic to
the trivial vector bundle over $\widetilde{A}^*\cong S^2\times
\R^3$
\[
\tau_{{\mathcal T}^{\widetilde{A}}\widetilde{A}^*}:{\mathcal
T}^{\widetilde{A}}\widetilde{A}^* \cong (S^2\times \R^3)\times
(\R^3\times \R^3)\to \widetilde{A}^*\cong S^2\times \R^3.
\]

Under the canonical identif\/ication $T\widetilde{A}^*\cong
TS^2\times (\R^3\times \R^3),$ the anchor map $\rho_{{\mathcal
T}^{\widetilde{A}} \widetilde{A}^*}:{\mathcal
T}^{\widetilde{A}}\widetilde{A}^*\to T\widetilde{A}^*$ of
${\mathcal T}^{\widetilde{A}}\widetilde{A}^*$ is given by
\[
\rho_{{\mathcal
T}^{\widetilde{A}}{\widetilde{A}^*}}((\vec{x},\pi),(\xi,\alpha))=((\vec{x},-\xi\times
\vec{x}),(\pi,\alpha))
\]
for $((\vec{x},\pi),(\xi,\alpha))\in{\mathcal
T}^{\widetilde{A}}\widetilde{A}^*\cong (S^2\times \R^3)\times
(\R^3\times \R^3).$ Moreover, we may choose, in a natural way, a
global basis $\{\widetilde{e}_i,\widetilde{f}_i\}_{i=1,2,3}$ of
$\Gamma({\mathcal T}^{\widetilde{A}}\widetilde{A}^*)$ such that
\[
\lcf\widetilde{e}_1,\widetilde{e}_2\rcf_{{\mathcal
T}^{\widetilde{A}}\widetilde{A}^*}=\widetilde{e_3}, \qquad
\lcf\widetilde{e}_3,\widetilde{e}_1\rcf_{{\mathcal
T}^{\widetilde{A}}\widetilde{A}^*}=\widetilde{e_2},\qquad
\lcf\widetilde{e}_2,\widetilde{e}_3\rcf_{{\mathcal
T}^{\widetilde{A}}\widetilde{A}^*}=\widetilde{e_1},
\]
and the rest of the fundamental Lie brackets are zero.

The symplectic section $\Omega_{{\mathcal
T}^{\widetilde{A}}\widetilde{A}^*}$ is given by
\[
\Omega_{{\mathcal
T}^{\widetilde{A}}\widetilde{A}^*}(\vec{x},\pi)((\xi,\alpha),(\xi',\alpha'))=\alpha'(\xi)-\alpha(\xi')
+\pi(\xi\times \xi'),
\]
for $(\vec{x},\pi)\in \widetilde{A}^*\cong S^2\times \R^3$ and
$(\xi,\alpha),(\xi',\alpha')\in\R^3\times \R^3.$

Now, we consider the Lie subgroup $G$ of $SO(3)$ of rotations about
the $z$-axis. The elements of $G$ are of the form
\[
A_\theta=\left (\begin{array}{ccc} \cos \theta&-\sin\theta &0\\
\sin\theta&\cos\theta&0\\0&0&1\end{array}\right) \qquad \mbox{with} \ \
\theta\in \R.
\]
It is clear that $G$ is isomorphic to $S^1=\R/(2\pi \Z).$

Note that the standard action of $G\cong S^1$ on $S^2$ is not free
(the north and the south poles are f\/ixed points). Therefore, we must
restrict this action to the open subset $M$ of $S^2$
\[
M=S^2-\{(0,0,1),(0,0,-1)\}.
\]
We have that the map $\mu:S^1\times \R\cong \R/(2\pi\Z)\times \R\to
M$ given by
\[
\mu([\theta],t)=\left(\frac{\cos \theta}{\cosh
t},\frac{\sin\theta}{\cosh t}, \tanh t\right),\qquad \mbox{for} \ \
([\theta],t)\in S^1\times \R
\]
is a dif\/feomorphism. Under this identif\/ication between $M$ and
$S^1\times \R,$ the action $\psi$ of $G\cong S^1$ on $M\cong
S^1\times \R$ is def\/ined by
\[
\psi([\theta'],([\theta],t))=([\theta+\theta'],t),
\qquad \mbox{for} \ \ [\theta']\in S^1 \ \ \mbox{and} \ \ ([\theta],t)\in S^1\times \R.
\]
Next, we consider the open subset of $\widetilde{A}$
\[
A=\tau_{\widetilde{A}}^{-1}(M)\cong (S^1\times \R)\times \R^3.
\]
It is clear that $A$ is an action Lie algebroid over $M\cong
S^1\times \R$ and the dual bundle to $A$ is the trivial vector
bundle $\tau_{A^*}:A^*\cong (S^1\times \R)\times \R^3\to S^1\times
\R.$ The restriction of $\widetilde{H}$ to $A^*$ is the
Hamiltonian function $H:A^*\cong (S^1\times \R)\times \R^3\to \R$
def\/ined by
\[
H(([\theta],t),\pi)=\frac{1}{2}\left(\frac{\pi_1^2}{I} +
\frac{\pi_2^2}{I}  + \frac{\pi_3^2}{J}\right) + (mgl) \tanh t.
\]

On the other hand, the action of $G\cong S^1$ on $M$ and the
standard action of $G$ on $\R^3$ induce an action $\Psi$ of $G$ on
$A$ in such a way $\Psi$ is a Lie algebroid action of $G$ on $A$.
Thus, we may consider the corresponding symplectic action
${\mathcal T}^*\Psi$ (see Example \ref{e3.6''}) of $G$ on the
symplectic Lie algebroid
\[
\tau_{{\mathcal T}^AA^*}:{\mathcal T}^AA^*\cong (M\times
\R^3)\times (\R^3\times \R^3)\to M\times \R^3.
\]
Note that $A^*$ is isomorphic to the open subset
$\tau_{\widetilde{A}^*}^{-1}(M)$ of $\widetilde{A}^*$ and that
\[
{\mathcal T}^AA^*\cong \tau^{-1}_{{\mathcal
T}^{\widetilde{A}}\widetilde{A}^*}(\tau_{\widetilde{A}^*}^{-1}(M)).
\]

So, the basis $\{\widetilde{e}_i,\widetilde{f}_i\}_{i=1,2,3}$ of
$\Gamma({\mathcal T}^{\widetilde A}\widetilde{A}^*)$ induces a
basis of $\Gamma({\mathcal T}^AA^*)$  which we will denote by
$\{\bar{e}_i,\bar{f}_i\}_{i=1,2,3}. $ The Hamiltonian section
${\mathcal  H}_H^{\Omega_{{\mathcal T}^{{A}}{A}^*}}$ of $H$ in the
symplectic Lie algebroid $({\mathcal T}^{A}{A}^*,$
$\Omega_{{\mathcal T}^{A}{A}^*})$ is given by
\begin{gather*}
{\mathcal H}_H^{\Omega_{{\mathcal
T}^{A}{A}^*}}(([\theta],t),\pi)=
\frac{{\pi}_1}{I}\bar{e}_1 +
\frac{\pi_2}{I}\bar{e_2} + \frac{\pi_3}{J}\bar{e}_3
+\left(\frac{(I-J)\pi_2\pi_3}{IJ}-\frac{mgl}{(\cosh t)}\sin\theta\right)\bar{f}_1\\
\phantom{{\mathcal H}_H^{\Omega_{{\mathcal
T}^{A}{A}^*}}(([\theta],t),\pi)=}{}-\left(\frac{(I-J)\pi_1\pi_3}{IJ}-\frac{mgl}{(\cosh
t)}\cos\theta\right)\bar{f}_2
\end{gather*}
for $(([\theta],t),\pi)\in A^*\cong (S^1\times \R)\times\R^3.$

Now, we consider the submanifold $N$ of $A^*\cong M\times \R^3$
\[
N=\{(([\theta],t),\pi)\in A^*/\pi_3=0\}
\]
and the Lie  subalgebroid $B$ (over $N$) of ${\mathcal T}^AA^*$
\[
B=\{((([\theta],t),\pi),(\xi,\alpha))\in {\mathcal
T}^AA^*/\pi_3=\alpha_3=0\}.
\]
If $\Omega_B$ is the restriction of $\Omega_{{\mathcal T}^AA^*}$
to the Lie subalgebroid $B$ then
\[
\dim (\ker\Omega_B(([\theta],t),\pi))=1,\;\;\;\; \mbox{ for all }
(([\theta],t),\pi)\in N
\]
and the section $s$ of $\tau_B:B\to N$ def\/ined by
\[
s=(\bar{e}_3 + \pi_2\bar{f}_1-\pi_1\bar{f}_2)_{|N}
\]
is a global basis of $\Gamma(\ker \Omega_B).$ Moreover, the
restriction of the symplectic action ${\mathcal T}^*\Psi$ to the
Lie algebroid $\tau_B:B\to N$ is a presymplectic action. In fact,
the action of $G$ on $N$ is given by
\[
([\theta'],(([\theta],t),\pi))\rightarrow
(([\theta+\theta'],t),A_{\theta'}\pi),
\]
for $[\theta']\in G\cong S^1$ and $(([\theta],t),\pi)\in N$, and
the action of $G$ on $B$ is
\[
([\theta'],((([\theta],t),\pi),(\xi,\alpha)))\to ((([\theta+\theta'],t),A_{\theta'}\pi),(A_{\theta'}\xi,A_{\theta'}\alpha)).
\]
Thus, using Theorem \ref{t2.3}, we may consider the reduced vector
bundle $\tau_{\widetilde{B}}:\widetilde{B}=(B/\ker \Omega_B)/G\to
\widetilde{N}=N/G.$ Furthermore, if $\widetilde{\pi}_B:B\to
\widetilde{B}$ and $\pi_N:N\to \widetilde{N}=N/G$ are the
canonical projections then a basis of the space
$\Gamma(B)_{\widetilde{\pi}_B}^p$ of
$\widetilde{\pi}_B$-projectable sections of $\tau_B:B\to
N$ is $\{e_1',e_2',s,f_1',f_2'\}$, where
\begin{gather*}
e_1'=(\cos\theta)\bar{e}_{1|N} + (\sin \theta)\bar{e}_{2|N},
\qquad e_2'=-(\sin\theta )\bar{e}_{1|N} +
\cos(\theta)\bar{e}_{2|N},
\\
f_1'=(\cos\theta)\bar{f}_{1|N} +  (\sin\theta)\bar{f}_{2|N},\qquad
f_2'=-(\sin\theta)\bar{f}_{1|N} +  (\cos\theta)\bar{f}_{2|N}.
\end{gather*}
We have that
\begin{gather*}
\lcf e_1',e_2'\rcf_B=(\sinh t)e_1' + s-
(\pi_2\cos\theta-\pi_1\sin\theta)f_1' + (\pi_1\cos\theta +
\pi_2\sin\theta)f_2',
\\
 \lcf e_1',f_1'\rcf_{B}=-(\sinh t)f_2',\qquad \lcf e_1',f_2'\rcf_B=(\sinh t)f_1',
\end{gather*}
and the rest of the Lie brackets between the elements of the basis
$\{e_1',e_2',s,f_1',f_2'\}$ are zero.

Note that the vertical bundle to $\pi_N$ is generated by the vector
f\/ield on $N$
\[
\frac{\partial}{\partial \theta}-\pi_2\frac{\partial}{\partial \pi_1} + \pi_1 \frac{\partial}{\partial \pi_2}
\]
and therefore, $\lcf e_1',e_2'\rcf_B, \;\lcf e_1', f_1'\rcf_B ,\;
\lcf e_1',f_2'\rcf_B\in \Gamma(B)_{\widetilde{\pi}_B}^p.$

Consequently, $\Gamma(B)^p_{\widetilde{\pi}_B}$ is a Lie
subalgebra of $(\Gamma(B),\lcf\cdot,\cdot\rcf_B)$  and
$\Gamma(\ker\Omega_B)$ is an ideal of
$\Gamma(B)^p_{\widetilde{\pi}_B}.$ This implies that the vector
bundle $\tau_{\widetilde{B}}:\widetilde{B}\to \widetilde{N}$
admits a unique Lie algebroid structure
$(\lcf\cdot,\cdot\rcf_{\widetilde{B}},\rho_{\widetilde{B}})$ such
that $\widetilde{\pi}_B$ is a Lie algebroid epimorphism. In
addition, $\Omega_B$ induces a symplectic section
$\Omega_{\widetilde{B}}$ on the Lie algebroid
$\tau_{\widetilde{B}}:\widetilde{B}\to \widetilde{N}$ such that
$\widetilde{\pi}_B^*\Omega_{\widetilde{B}}=\Omega_B$ (see Theorem
\ref{t2.4}).

On the other hand, using Theorem \ref{t2.6}, we deduce that the
original Hamiltonian system, with Hamiltonian function $H:A^*\to
\R$ in the symplectic Lie algebroid $\tau_{{\mathcal
T}^{A}A^*}:{\mathcal T}^{A}A^*\to A^*$, may be reduced to a
Hamiltonian system, with Hamiltonian function
$\widetilde{H}:\widetilde{N}\to \R$, in the symplectic Lie
algebroid $\tau_{\widetilde{B}}:\widetilde{B} \to \widetilde{N}$.

Next, we will give an explicit description of the reduced symplectic
Lie algebroid and the reduced Hamiltonian system on it.

In fact, one may prove that the Lie algebroid
$\tau_{\widetilde{B}}:\widetilde{B}\to \widetilde{N}$ is
isomorphic to the trivial vector bundle $\tau_{\R^3\times
\R^4}:{\R^3\times \R^4}\to \R^3$ with basis $\R^3$ and f\/iber
$\R^4$ and, under this identif\/ication, the anchor map
$\rho_{\widetilde{B}}$ is given by
\begin{gather*}
\rho_{\widetilde{B}}((t,v),(\eta,\beta))=-\eta_2(\cosh
t) \frac{\partial}{\partial t} + (\beta_1-\eta_1\nu_2
\sinh t) \frac{\partial} {\partial \nu
_1} +(\beta_2 + \eta_1\nu_1 \sinh
t) \frac{\partial}{\partial \nu_2},
\end{gather*}
for $(t,\nu \equiv (\nu_1,\nu_2))\in \R\times \R^2\cong \R^3$ and
$(\eta\equiv(\eta_1,\eta_2),\beta\equiv (\beta_1,\beta_2))\in
\R^2\times \R^2\cong \R^4.$

Moreover, if $\{e_i\}_{i=1,\dots, 4}$ is the canonical basis of
$\Gamma(\widetilde{B})\cong \Gamma(\R^3\times \R^4)$ we have that
\begin{gather*}
\lcf e_1,e_2\rcf_{\widetilde{B}}=(\sinh t)e_1 -\nu_2e_3 +
\nu_1e_4,
\\
\lcf e_1,e_3\rcf_{\widetilde{B}}=-(\sinh t)e_4,\qquad \lcf
e_1,e_4\rcf_{\widetilde{B}}=(\sinh t)e_3
\end{gather*}
and the rest of the Lie brackets between the elements of the basis
is zero.

Finally, one may prove that the reduced symplectic section
$\Omega_{\widetilde{B}}$, the reduced Hamiltonian function
$\widetilde{H}$ and the corresponding reduced Hamiltonian vector
f\/ield ${\mathcal H}_{{H}}^{\{\cdot,\cdot\}_{\widetilde{N}}}$ on
$\widetilde{N}\cong \R^3$ are given by the following expressions
\begin{gather*}
\Omega_{\widetilde{B}}(t,\nu)((\eta,\beta),(\eta',\beta'))=\beta'(\eta)-\beta(\eta'), \\
\widetilde{H}(t,\nu)= \frac{1}{2} \left(\frac{\nu_1^2}{I}
+ \frac{\nu_2^2}{I}\right) + mgl(\tanh t),\\
{\mathcal H}_{\widetilde{H}}^{\{\cdot,\cdot\}_{\widetilde{N}}}=
- \frac{\nu_2}{I} \cosh t \frac{\partial
}{\partial t} -  \frac{\nu_1\nu_2}{I} \sinh t
 \frac{\partial }{\partial \nu_1}
  + \left( \frac{mgl}{\cosh t} +
\frac{\nu_1^2}{I} \sinh  t\right)
\frac{\partial }{\partial \nu_2}.
\end{gather*}

\section{Conclusions and future work}

In this paper we have generalized the Cartan symplectic reduction
by a submanifold in the presence of a symmetry Lie group to the
Lie algebroid setting. More precisely, we develop a~procedure to
reduce a symplectic Lie algebroid to a certain quotient symplectic
Lie algebroid, using a Lie subalgebroid and a symmetry Lie group.
In addition, under some mild assumptions we are able to reduce the
Hamiltonian dynamics. Several examples illustrate our theory.

As we already mentioned in the introduction, a particular example
of reduction of a symplectic manifold is the well known
Marsden--Weinstein reduction of a symmetric Hamiltonian dynamical
system $(M,\Omega ,H)$ in the presence of a $G$-equivariant
momentum map $J:M\to \mathfrak g^*$, where $\mathfrak g$ is the
Lie algebra of the Lie group $G$ which acts on $M$ by
$\Phi:G\times M\to M$. In this situation, the submanifold is the
inverse image $J^{-1}(\mu)$ of a regular value $\mu \in\mathfrak
g^*$ of $J$, and the Lie group acting on $J^{-1}(\mu)$ is the
isotropy group of $\mu$. It would be interesting to generalize
this procedure to the setting of symplectic Lie algebroids. In
this case, we would have to def\/ine a proper notion of a momentum
map for this framework. This appropiate moment map would allow us
to get the right submanifold in order to apply the reduction
procedure we have developed in this paper. Moreover, under these
hypotheses, we should be able to reduce also the Hamiltonian
dynamics. These topics are the subject of a forthcoming paper (see~\cite{IMMMP}).

A natural generalization of symplectic manifolds is that of
Poisson manifolds. This type of manifolds play also a prominent
role in the Hamiltonian description of Mechanics, in particular
for systems with constraints or with symmetry. For Poisson
manifolds, one can also develop a~reduction procedure (see~\cite{MR}). On the other hand, as we mentioned in Section~2.2, a
symplectic section on a Lie algebroid $A$ can be seen as a
particular example of a triangular matrix for~$A$, that is, a
Poisson structure in the Lie algebroid setting. Therefore, it
should be interes\-ting to generalize our reduction constructions to
Poisson structures on Lie algebroids and, more generally, to
obtain a reduction procedure for Lie bialgebroids (we remark that
triangular matrices on Lie algebroids are a particular example of
Lie bialgebroids, the so-called triangular Lie bialgebroids).

On the other hand, recently, in \cite{BCG} it has been described a
reduction procedure for Courant algebroids and Dirac structures.
We recall that there is a relation between Courant algebroids and
Lie bialgebroids: there is a one-to-one correspondence between Lie
bialgebroids and pairs of complementary Dirac structures for a
Courant algebroid (see \cite{LWX}). Since symplectic Lie
algebroids are a particular example of this situation, it would be
interesting to compare both approaches. This is the subject of a
forthcoming paper.

\pdfbookmark[1]{Appendix}{Appendix}
\section*{Appendix}

\newcounter{apendice}
\def\theequation{\Alph{apendice}.\arabic{equation}}
\def\thetheorem{\Alph{apendice}.\arabic{theorem}}
\def\theremark{\Alph{apendice}.\arabic{remark}}

\setcounter{equation}{0} \setcounter{apendice}{1}
\setcounter{theorem}{0} \setcounter{remark}{1}

Suppose that $\tau_A:A\to M$ and
$\tau_{\widetilde{A}}:\widetilde{A}\to \widetilde{M}$ are real
vector bundles and that $\pi_A:A\to \widetilde{A}$ is a~vector
bundle epimorphism over the surjective submersion $\pi_M:M\to
\widetilde{M}.$

Thus, one may consider the vector subbundle $\ker\pi_A$ of $A$
whose f\/iber over the point $x\in M$ is the kernel of the linear
epimorphism $(\pi_A)_{|A_x}:A_x\to \widetilde{A}_{\pi_M(x)}.$

On the other hand, a section $X$ of the vector bundle $\tau_A:A\to
M$ is said to be $\pi_A$-projectable if there exists a section
$\pi_A(X)$ of $\tau_{\widetilde{A}}:\widetilde{A}\to
\widetilde{M}$ such that
\begin{equation}\label{secproj}
\pi_A\circ X=\pi_A(X)\circ \pi_M.
\end{equation}
We will denote by $\Gamma(A)_{\pi_A}^p$ the set of
$\pi_A$-projectable sections. $\Gamma(A)_{\pi_A}^p$ is a
$C^\infty(\widetilde{M})$-module and it is clear that
$\Gamma(\ker\pi_A)$ is a $C^\infty(\widetilde{M})$-submodule of
$\Gamma(A)^p_{\pi_A}.$

Note that if $X\in \Gamma(A)^p_{\pi_A}$ then there exists a unique
section $\pi_A(X)$ of the vector bundle such that (\ref{secproj})
holds.

Thus, we may def\/ine a map
\[
\Pi_A:\Gamma(A)^p_{\pi_A}\to \Gamma(\widetilde{A})
\]
and it follows that $\Pi_A$ is an epimorphism of
$C^\infty(\widetilde{M})$-modules and that
$\ker\Pi_A=\Gamma(\ker\pi_A)$. Therefore, we have that the
$C^\infty(\widetilde{M})$-modules $\Gamma(\widetilde{A})$ and
$\displaystyle\frac{\Gamma(A)^p_{\pi_A}}{\Gamma(\ker\pi_A)}$ are
isomorphic, that is,
\begin{equation}\label{A1}
\Gamma(\widetilde{A})\cong
\frac{\Gamma(A)^p_{\pi_A}}{\Gamma(\ker\pi_A)}.
\end{equation}
Moreover, one may prove that following result.

\begin{theorem}\label{tA1}
Let $\pi_A:A\to \widetilde{A}$ be a vector bundle epimorphism
between the vector bundles $\tau_A:A\to M$ and
$\tau_{\widetilde{A}}:\widetilde{A}\to \widetilde{M}$ over the
surjective submersion $\pi_M:M\to \widetilde{M}$ and suppose that
$A$ is a Lie algebroid with Lie algebroid structure
$(\lcf\cdot,\cdot\rcf_A,\rho_A).$ Then, there exists a unique Lie
algebroid structure on the vector bundle
$\tau_{\widetilde{A}}:\widetilde{A}\to \widetilde{M}$ such that
$\pi_A$ is a Lie algebroid epimorphism if and only if the
following conditions hold:

\begin{enumerate}\itemsep=0pt
\item[i)] The space $\Gamma(A)^p_{\pi_A}$ is a Lie subalgebra of
the Lie algebra $(\Gamma(A),\lcf\cdot,\cdot\rcf_A)$ and \item[ii)]
$\Gamma(\ker \pi_A)$ is an ideal of the Lie algebra
$\Gamma(A)^p_{\pi_A}.$
\end{enumerate}
\end{theorem}
\begin{proof} Assume that
$(\lcf\cdot,\cdot\rcf_{\widetilde{A}},\rho_{\widetilde{A}})$ is a
Lie algebroid structure on the vector bundle
$\tau_{\widetilde{A}}:\widetilde{A}\to \widetilde{M}$ such that
$\pi_A$ is a Lie algebroid epimorphism.

If $X\in \Gamma(A)^p_{\pi_A}$ then, using that
\[
d^A(\widetilde{f}\circ
\pi_M)=(\pi_A)^*(d^{\widetilde{A}}\widetilde{f}),\qquad \forall\,
\widetilde{f}\in C^\infty(\widetilde{M}),
\]
we directly conclude that the vector f\/ield $\rho_A(X)$ on $M$ is
$\pi_M$-projectable on the vector f\/ield
$\rho_{\widetilde{A}}(\pi_A(X))$ on $\widetilde{M}$.

Thus, if $Y\in \Gamma(A)^p_{\pi_A}$ and $\widetilde{\alpha}\in
\Gamma(\widetilde{A}^*)$, we have that
\begin{gather*}
\widetilde{\alpha}(\pi_M(x))(\lcf\pi_A(X),\pi_A(Y)\rcf_{\widetilde{A}}(\pi_M(x)))\\
\qquad{}=-(\pi_A^*(d^{\widetilde{A}}
\widetilde{\alpha}))(x)(X(x),Y(x))+\rho_A(X)(x)(\pi_A^*\widetilde{\alpha})(Y))
-\rho_A(Y)(x)((\pi_A^*\widetilde{\alpha})(X)),
\end{gather*}
for all $x\in M$. Therefore, since
$(\pi_A^*(d^{\widetilde{A}}\widetilde{\alpha}))=
d^A(\pi_A^*\widetilde{\alpha})$ it follows that
\begin{gather*}
\widetilde{\alpha}(\pi_M(x))(\lcf \pi_A(X),
\pi_A(Y)\rcf_{\widetilde{A}}(\pi_M(x)))\\
\qquad {}=(\pi_A^*\widetilde{\alpha})(x)(\lcf X,Y\rcf_A(x))=
\widetilde{\alpha}(\pi_M(x))(\pi_A \lcf X,Y\rcf_A(x)).
\end{gather*}
Consequently, we deduce the following result
\[
X,Y\in \Gamma(A)^p_{\pi_A}\Rightarrow \pi_A\circ \lcf
X,Y\rcf_A=\lcf\pi_A(X), \pi_A(Y)\rcf_{\widetilde{A}}\circ \pi_M.
\]
This proves $(i)$ and $(ii)$.

Conversely, suppose that conditions $(i)$ and $(ii)$ hold. Then,
we will def\/ine a Lie bracket $\lcf\cdot,\cdot\rcf_{\widetilde{A}}$
on the space $\Gamma(\widetilde{A})$ as follows. If
$\widetilde{X}$ and $\widetilde{Y}$ are sections of
$\tau_{\widetilde{A}}:\widetilde{A}\to \widetilde{M}$, we can
choose $X,Y\in \Gamma(A)^p_{\pi_A}$ such that
\[
\pi_A(X)=\widetilde{X},\qquad \pi_A(Y)=\widetilde{Y}.
\]
Now, we have that $\lcf X,Y\rcf_A\in\Gamma(A)^p_{\pi_A}$ and we
def\/ine
\begin{equation}\label{Bractil}
\lcf\widetilde{X},\widetilde{Y}\rcf_{\widetilde{A}}=\lcf\pi_A(X),
\pi_A(Y)\rcf_{\widetilde{A}}=\pi_A \lcf X,Y\rcf_A.
\end{equation}

Using condition $(ii)$, it follows that the map
$\lcf\cdot,\cdot\rcf_{\widetilde{A}}: \Gamma(\widetilde{A})\times
\Gamma(\widetilde{A}) \to \Gamma(\widetilde{A})$ is well-def\/ined.
Moreover, since $\lcf\cdot,\cdot\rcf_A$ is a Lie bracket on
$\Gamma(A)$, we conclude that
$\lcf\cdot,\cdot\rcf_{\widetilde{A}}$ is also a Lie bracket on~$\Gamma(\widetilde{A}).$

Next, will see that
\begin{equation}\label{Anch}
Z\in \Gamma(\ker\pi_A) \ \Rightarrow \ \rho_A(Z) \mbox{ is a
$\pi_M$-vertical vector f\/ield.}
\end{equation}
In fact, if $\widetilde{f}\in C^\infty(\widetilde{M})$ we have that
\[
\lcf Z,(\widetilde{f}\circ \pi_M)Y\rcf_A\in
\Gamma(\ker\pi_A),\qquad \mbox{for all} \ \ Y\in
\Gamma(A)^p_{\pi_A}.
\]
On the other hand,
\[
\lcf Z,(\widetilde{f}\circ \pi_M)Y\rcf_A=(\widetilde{f}\circ
\pi_M)\lcf Z,Y\rcf_A + \rho_A(Z)(\widetilde{f}\circ \pi_M)Y
\]
which implies that $\rho_A(Z)(\widetilde{f}\circ \pi_M)=0$ and,
thus, the vector f\/ield $\rho_A(Z)$ is vertical with respect to the
map $\pi_M$.

Now, let $\widetilde{X}$ be a section of
$\tau_{\widetilde{A}}:\widetilde{A}\to \widetilde{M}$ and $X$ be a
section of $\tau_A:A\to M$ such that $\pi_A(X)=\widetilde{X}.$ If
$\widetilde{f}\in C^\infty(\widetilde{M})$ it follows that
\[
\lcf X, (\widetilde{f}\circ \pi_M)Y\rcf_A\in
\Gamma(A)^p_{\pi_A},\qquad \mbox{for all} \ \ Y \in
\Gamma(A)^p_{\pi_A}
\]
and, since
\[
\lcf X,(\widetilde{f}\circ \pi)Y\rcf_A=(\widetilde{f}\circ
\pi_M)\lcf X,Y\rcf_A + \rho_A(X)(\widetilde{f}\circ \pi_M)Y
\]
we obtain that $\rho_A(X)(\widetilde{f}\circ \pi_M)$  is a basic
function with respect to the map $\pi_M$. Therefore, the vector
f\/ield $\rho_A(X)$ is $\pi_M$-projectable to a vector f\/ield
$\rho_{\widetilde{A}}(\widetilde{X})$ on $\widetilde{M}.$

Using (\ref{Anch}), we deduce that the map $
\rho_{\widetilde{A}}:\Gamma(\widetilde{A})\to {\mathfrak
X}(\widetilde{M})$ is well-def\/ined. In fact,
$\rho_{\widetilde{A}}$ is a~homomorphism of
$C^\infty(\widetilde{M})$-modules and if we also denote by
$\rho_{\widetilde{A}}:\widetilde{A}\to T\widetilde{M}$ the
corresponding bundle map then
\begin{equation}\label{Anchtil}
T\pi_M\circ \rho_A=\rho_{\widetilde{A}}\circ \pi_A.
\end{equation}
In addition, since $(\lcf\cdot,\cdot\rcf_A,\rho_A)$ is a Lie
algebroid structure on $A$, we conclude that the pair
$(\lcf\cdot,\cdot\rcf_{\widetilde{A}}, \rho_{\widetilde{A}})$ is
also a Lie algebroid structure on $\widetilde{A}$.

On the other hand, from (\ref{Anchtil}), it follows that
\[
\pi_A^*(d^{\widetilde{A}}\widetilde{f})=d^A(\widetilde{f}\circ
\pi_M),\qquad \mbox{for} \ \ \widetilde{f}\in
C^\infty(\widetilde{M}).
\]
Furthermore, using (\ref{Bractil}) and (\ref{Anchtil}), we obtain that
\[
\pi_A^*(d^A\widetilde{\alpha})=d^A(\pi_A^*\widetilde{\alpha}),\qquad \mbox{for} \ \ \widetilde{\alpha}\in \Gamma(\widetilde{A}^*).
\]
Consequently, $\pi_A$ is a Lie algebroid epimorphism.

Now, suppose that
$(\lcf\cdot,\cdot\rcf_{\widetilde{A}}',\rho_{\widetilde{A}}')$ is
another Lie algebroid structure on the vector bundle
$\tau_{\widetilde{A}}:\widetilde{A}\to \widetilde{M}$ such that
$\pi_A$ is a Lie algebroid epimorphism. Denote by
$d^{'\widetilde{A}}$ the dif\/ferential of the Lie algebroid
$(\widetilde{A},\lcf\cdot,\cdot\rcf'_{\widetilde{A}},\rho'_{\widetilde{A}})$.
Then, the condition
\[
\pi_A^*(d^{'\widetilde{A}}\widetilde{f})= d^A(\widetilde{f}\circ
\pi_M),\qquad \mbox{for} \quad \widetilde{f}\in
C^\infty(\widetilde{M}),
\]
implies that $T\pi_M\circ \rho_A=\rho'_{\widetilde{A}}\circ
\pi_A.$ Thus, from (\ref{Anchtil}), we deduce that
$\rho_{\widetilde{A}}= \rho'_{\widetilde{A}}.$

Therefore, the condition
\[
\pi_A^*(d^{'\widetilde{A}}\widetilde{\alpha})=
d^A(\pi_A^*\widetilde{\alpha}),\qquad \mbox{for} \ \
\widetilde{\alpha} \in \Gamma(\widetilde{A}^*),
\]
implies that
\[
\lcf \pi_A(X),\pi_A(Y)\rcf_{\widetilde{A}}' =\pi_A(\lcf
X,Y\rcf_A),\qquad \mbox{for} \ \ X,Y\in \Gamma(A)^p_{\pi_A}
\]
and, using (\ref{Bractil}), we conclude that
$\lcf\cdot,\cdot\rcf_{\widetilde{A}}=\lcf\cdot,\cdot\rcf_{\widetilde{A}}'$.
\end{proof}

\begin{remark}
An equivalent dual version of Theorem~\ref{tA1} was proved in
\cite{CNS}.
\end{remark}

\subsection*{Acknowledgements}
Supported in part by BFM2003-01319, MTM2006-03322, MTM2004-7832,
S-0505/ESP/0158 of the CAM and BFM2003-02532. J.C.~Marrero wishes to
thank the organizers of the workshop Geometric Aspects of
Integrable Systems (Coimbra, 2006) for their hospitality and their
job in the organization. Finally, we also would like to thank the
anonymous referees for the observations that have improved the
manuscript.

\pdfbookmark[1]{References}{ref}
\LastPageEnding
\end{document}